\documentclass[11 pt]{amsart}
\usepackage[dvips]{graphicx}
\usepackage{amssymb}
\usepackage[english]{babel}
\usepackage{amsfonts,amsthm,amsmath}
\usepackage{epsfig}
\input{amssym.def}

\newtheorem{theorem}{Theorem}
\newtheorem{prop}[theorem]{Proposition}
\newtheorem{proposition}[theorem]{Proposition}
\newtheorem{lemma}[theorem]{Lemma}
\newtheorem{corollary}[theorem]{Corollary}

\newtheorem{remark}[theorem]{Remark}

\newcommand{\RR}{\mathbb{R}}

\begin{document}

\title[Inverse spectral problems.]{Inverse spectral problems for Schr\"odinger and pseudo-differential
operators.}
\author{Brice Camus}
\address{Ludwig Maximilians Universit\"at M\"unchen,\\
Mathematisches Institut, Theresienstr. 39 D-80803
M\"unchen.\\Email: camus@math.lmu.de} \maketitle
\begin{abstract}
\noindent Starting from the semi-classical spectrum of
Schr\"odinger operators $-h^2\Delta+V$ (on $\mathbb{R}^n$ or on a
Riemannian manifold) it is possible to detect critical levels of
the potential $V$. Via micro-local methods one can express
spectral statistics in terms of different invariants:
\begin{itemize}
\item Geometry of energy surfaces (heat invariant like).
\item Classical orbits (wave invariants).
\item But also classical equilibria (new wave invariants).
\end{itemize}
Any critical point of $V$ with zero momentum is an equilibrium of
the flow and generates many singularities in the semi-classical
distribution of eigenvalues. Via sharp spectral estimates, this
phenomena indicates the presence of a critical energy level and
the information contained in this singularity allows to
reconstruct partially the local shape of $V$. Several generalizations of this approach are also proposed.\medskip
\end{abstract}
\section{Introduction.}
\subsection{Background and basic definitions.}
In this article we are here interested in the inverse spectral
problem for partial differential operators and pseudo-differential
operators in the semi-classical or high-energy regime. A natural
question is to try to understand how the semi-classical spectrum
of such an operator can describe the shape of the graph of the
(principal) symbol: critical points, extrema and associated local
Taylor expansions. Because the spectrum is invariant under
translation of the symbol it is in general not possible to obtain
more than a qualitative answer. For example there is no hope to
locate critical points of the symbol starting only form the
spectrum.

The results we would like to present are perhaps of particular
interest for $h$-quantized Schr\"odinger operators on
$L^2(\mathbb{R}^n)$:
\begin{equation*}
P_{h}=-h^2\Delta+V,
\end{equation*}
also called semi-classical Schr\"odinger operators. Here we will
assume that the potentials $V$ are smooth on $\mathbb{R}^n$ and
bounded from below. For this class of operators, the question is
then to understand how certain fluctuations in the semi-classical
spectrum can describe the shape of $V$. But our results will also
apply to more general operators like $h$-pseudo-differential,
$h$-admissible operators (see definitions below) or Schr\"odinger
operators on a Riemannian manifold $M$ (simply replace the Laplace
operator by a Laplace-Beltrami operator $\Delta_M$, see section
5). Most of these modifications are possible because our methods
are micro-local and do not
use global results on these operators.\medskip\\
\textbf{Notations.} Before entering into the details, let us give
some definitions and recall some basic facts about the spectral
theory for $P_h$. By a standard result, see \cite{Ber-Shu}, when
$V$ is bounded from below and with tempered growth, $P_{h}$ has a
self-adjoint realization on a dense subset of $L^2(\mathbb{R}^n)$.
As usually, to this quantum operator $P_{h}$ we can associate a
classical counterpart with the Hamiltonian function
$p(x,\xi)=\xi^2+V(x)$, or total energy, on the phase space
$\mathbb{R}^n\times\mathbb{R}^n$. In what follows, we note
$\Phi_{t}$ the flow  of the Hamiltonian vector field:
\begin{equation*}
H_{p}=\partial _{\xi}p.\partial_{x} -\partial _{x}p.\partial_{\xi}.
\end{equation*}
A classical energy surface is:
\begin{equation*}
\Sigma_E=\{ (x,\xi)\in \mathbb{R}^n\times\mathbb{R}^n \text{ / }
\xi^2+V(x)=E\}
\end{equation*}
and for a classical energy $E\in\mathbb{R}$ we will say that:
\begin{itemize}
\item $E$ is regular when $dp\neq 0$ everywhere on $\Sigma_E$.
\item $E=E_c$ is critical if $dp=0$ somewhere on $\Sigma_{E_c}$.
\end{itemize}
In the present article we are very precisely interested in a
relation between the asymptotic properties as $h\rightarrow 0^+$
of eigenvalues $\lambda_j(h)$ of $P_h$ :
\begin{equation*}
P_{h}\varphi_j(x,h)=\lambda_j(h) \varphi_j(x,h), \text{ }
\varphi_j\in L^2(\mathbb{R}^n),
\end{equation*}
and the set of fixed point for the map $\Phi_{t}$ (viewed as a map
on $\mathbb{R}_t\times T^*\mathbb{R}^n$). We recall in section 2,
a sufficient condition to get discrete spectrum. It is well known,
see \cite{Cdv}, that this semi-classical problem and the
high-energy limit for the spectrum of the Laplace operator
$\Delta_M>0$ on a compact Riemannian-manifold $M$ are related.
This relation can be viewed by quantizing $1/h=\sqrt{\lambda_j}$
where $\lambda_j\rightarrow +\infty$ is the increasing sequence of
eigenvalues of $\Delta_M$. But in this article we
will mainly consider the semi-classical problem and associated micro-local methods.\medskip\\
\textbf{Duality between the quantum and classical worlds.} In
geometry spectrum and periodic orbits can be related, in a very
explicit way, by means of the Selberg (e.g. see \cite{Hej}) and
Duistermaat-Guillemin \cite{D-G} trace formulae. It is in general
in the most atypical situations, compact surfaces of constant
negative curvatures or at the opposite completely integrable
systems (e.g. like the free Laplacian on a flat torus) that the
most explicit and exact results can be obtained. In quantum
mechanics, the existence of such a relation is strongly suggested
by the correspondence principle which asserts that, in the
semiclassical regime $h \rightarrow 0$, many properties of $P_{h}$
can be related to integral curves of $\Phi_t$ and many invariant
attached to the flow around these curves (see
Eq.(\ref{Gutzwiller}) below). For a general Hamiltonian $H$, not
necessarily of the form kinetic energy plus potential, this
correspondence principle is also true under very reasonable
assumptions on the symbol $h(x,\xi)$ of $H$, a function on the
phase space.

In physics a more precise formulation of this principle appeared
in the works of Balian\&Bloch \cite{BB} and Gutzwiller \cite{GUT}.
The Gutzwiller formula is usually written as a trace formula for
the resolvent of $P_h$ at a given energy $E$:
\begin{equation} \label{Gutzwiller}
\sum\limits_{j\in\mathbb{N}} \frac{1}{\lambda_j(h)-E}=
\frac{\mathrm{Lvol}(\Sigma_E)} {(2\pi h)^n} +\frac{1}{ih}
\sum\limits_{\gamma \in \Sigma_E} A_\gamma
e^{\frac{i}{h}S_\gamma},
\end{equation}
where in the r.h.s the sum concerns the closed orbits $\gamma$
inside $\Sigma_E$. Here $\mathrm{Lvol}(\Sigma_E)$ is the
Riemannian volume of $\Sigma_E$ (defined w.r.t. the invariant
Liouville measure), $S_\gamma=\int_\gamma \xi dx$ and $A_\gamma$
are respectively the classical action and the stability factor
(including the Maslov phase) of the curve $\gamma$. Recall that
the Liouville-volume $\mathrm{Lvol}(\Sigma_E)$ satisfies the
co-area formula:
\begin{equation*}
\int\limits_{[a,b]} \mathrm{Lvol}(\Sigma_E) dE
=\int\limits_{\{(x,\xi): a<p(x,\xi)<b\}}
dxd\xi=\mathrm{Vol}_{\mathbb{R}^{2n}} (p^{-1}([a,b]),
\end{equation*}
where the measure on the r.h.s. is the Lebesgue measure of the
pull-back.

In mathematics and in physics, such a relation between spectrum
and periodic orbits provides a powerful tool of analysis and
computation. See e.g. \cite{Laz} concerning the asymptotic
behavior of eigenvectors $\varphi_j(x,h)$ and \cite{Haa} for
various applications in quantum chaos. See also \cite{Cdv} for a nice overview and applications in Riemannian geometry.\medskip\\
\textbf{Mathematical problems.} For a Schr\"odinger operator on $\mathbb{R}^n$, it is easy to check that two different type of divergence generally occur in
Eq.(\ref{Gutzwiller}) :\medskip\\
\textbf{1)} The sum over the spectrum is divergent when the
resolvent is not a trace-class operator. In particular this is the
case if $V$ does not go fast enough to $\infty$ when
$|x|\rightarrow \infty$. When the sum appears to be convergent it
can also have a divergent behavior when $h\rightarrow 0$.
Worst, it can be that both sides of Eq. (\ref{Gutzwiller}) do not fit in the regime $h\rightarrow 0$.\\
\textbf{2)} The sum over closed orbits is generally divergent.
This is the case if $|A_\gamma|$ does not decrease fast enough or
if the number of periodic orbits of period smaller than $T$ is exponentially growing with $T$.\medskip\\
For example, if $n=1$, using scaling and the asymptotic properties
of the spectrum (here simply given by some Bohr-Sommerfeld
quantization conditions, see \cite{Ber-Shu}) it is easy to check
that the trace of the resolvent of :
\begin{equation*}
Q_h=-h^2\frac{d^2}{dx^2} +|x|^\alpha, \, h>0,\, \alpha>0,
\end{equation*}
exists if and only if $\alpha>2$. Here, the harmonic oscillator,
obtained for $\alpha=2$ for which $\lambda_j(h)=h(2j+1)$, is the
limit case and the series diverges like the harmonic series.
\subsection{Mathematical approach of the Gutzwiller
formula.}$\,$\\
As seen above, the question to remove divergences has many
important implications\footnote{Many important questions
concerning the range of trace formulae (e.g., their validity
beyond the Ehrenfest-time) are still open. We do not discuss these
questions in this article.} and we explain now a mathematical way
to solve this problem via a smoothing of the so-called spectral
density. To begin the discussion, simply assume that:
\begin{quote}
For some $E\in \mathbb{R}$, the spectrum of $P_h$ is discrete,
with finite multiplicities, in the interval
$[E-\varepsilon,E+\varepsilon ]$, $\varepsilon>0$.
\end{quote}
A sufficient condition to obtain this property, uniformly w.r.t.
$E$, is given in section 2. A well-posed problem is to study the
asymptotic behavior of the spectral distributions:
\begin{equation}
\Upsilon (E,h,\varphi )=\sum\limits_{|\lambda _{j}(h)-E|\leq
\varepsilon }\varphi (\frac{\lambda _{j}(h)-E}{h}), \text{
as } h\rightarrow 0, \label{Def trace}
\end{equation}
where $\varphi$ is a test function chosen to remove the
divergences. We can justify this terminology if we observe that
the truncated spectral distribution:
\begin{equation*}
T_{E,\varepsilon}(x)=\sum\limits_{|\lambda_j(h)-E|\leq
\varepsilon} \delta_{\lambda_j(h)}(x),\, \langle \delta_{x_0},
f\rangle = f(x_0),
\end{equation*}
acting on a function $\varphi$ shifted by $E$ and scaled w.r.t.
$h$ provides :
\begin{equation*}
\Upsilon(E,h,\varphi)=\left\langle T_{E,\varepsilon}(x) ,
\varphi(\frac{x-E}{h}) \right \rangle.
\end{equation*}
In reality this scaling w.r.t. $h$ is very important and is used
to get parametrices involving the classical dynamics. Also it is
not very hard to verify that when $\varphi\in\mathcal{S}(\RR)$ the
size of the truncation (materialized here as $\varepsilon$) is
irrelevant on a scale of size $\mathcal{O}(h^\infty)$ as long as
$\varepsilon$ is strictly positive. I refer to section 3 for these
points but I simply recall that $\mathcal{O}(h^\infty)$ is the
class of functions of $h$ being in
$\mathcal{O}(h^k)$ for every  $k\in\mathbb{N}$ near $h=0$.\medskip\\
\textbf{Statistical quantum mechanics.} In general, apart in some
very specific situations, it is not possible to compute explicitly
the spectrum of $P_h$ and a motivation to do semi-classical or
high-energy estimates is to derive statistics about eigenvalues
and their distribution. For example, in Eq.(\ref{Def trace}) the
formal choice of $\varphi$ as the characteristic function of
$[-\eta,\eta]$, $0<\eta<\varepsilon$, determines the number
$N(E,h)$ of bound states in $[E-\eta h,E+\eta h]$.

This formal correspondence between $\Upsilon$ and the micro-local
counting function $N$ has a mathematically rigorous formulation in
term of Tauberian-theorems, see e.g. \cite{BPU}. Under certain
(generic) conditions\footnote{In particular the condition that $E$
is non-critical for $p$, see below.} on the symbol $p$ of $P_h$,
it can be proven that $N(E,h)$ is proportional to $h^{1-n}$ times
the Liouville-volume of the energy shell $\Sigma_E$:
\begin{equation*}
N(E,h)\sim \frac{1}{(2\pi h)^{n-1}} \mathrm{Lvol}(\Sigma_E)\text{
as } h\rightarrow 0^+.
\end{equation*}
This is a micro-local formulation of the Weyl-law. A fortiori when
$n\geq 2$ this implies that the finite sum defining
$\Upsilon(E,h,\varphi)$ will involve a large number of eigenvalues
in the regime $h\rightarrow 0^+$. In general the formula for
$N(E,h)$ can be formally integrated to obtain a formula for
counting eigenvalues in a compact interval:
\begin{equation*}
\mathcal{N}([a,b],h)=\# \{ j\in \mathbb{N}: \lambda_j(h)\in
[a,b]\} \sim \frac{1}{(2\pi h)^h} \mathrm{Vol}_{\RR^{2n}}
p^{-1}([a,b]).
\end{equation*}
Of course if the operator is bounded from below you can also use
$\mathcal{N}(x,h):= \mathcal{N}(]-\infty,x],h)$. The question to
estimate the remainder function for $\mathcal{N}$ is in general a
relatively complicated problem and requires to use the properties
of the underlying classical dynamics inside $p^{-1}([a,b])$.
Several other related problems, Riesz-moments or Lieb-Thirring
inequalities, can be formulated in terms of $N(E,h)$ and these
problems are important in the 'stability of matter' problem. See
\cite{Lieb} for an overview and references. Finally, I mention
that certain Schr\"odinger operators with very singular critical
sets (e.g. see \cite{BPU,Cam3}) or non-confining potentials (e.g.
see \cite{Sim2}) can lead to some very different kind of
'Weyl-asymptotics' for $N$ or $\mathcal{N}$.
\medskip\\
\textbf{Relation with the classical dynamics.} In reality the
quantity $\Upsilon(E,h,\varphi)$ contains many interesting
information (a priori more than the counting functions) since the
asymptotic expansion of $\Upsilon(E,h,\varphi)$ when $h\rightarrow
0^+$ involves explicitly the classical dynamics on $\Sigma_E$ and
in particular the set of fixed point of the flow inside the energy
surface:
\begin{equation*}
\mathrm{Fix}_E=\{ (T,x,\xi) \in \RR\times \Sigma_E :
\Phi_T(x,\xi)=(x,\xi)\}.
\end{equation*}
We recall that $E$ is regular if $\nabla p(x,\xi )\neq 0$ on
$\Sigma_{E}$ and critical otherwise. Every critical point
$(x_0,\xi_0)\in\Sigma_{E_c}$ of $p$ is a fixed point of our flow
$\Phi_t$ since $H_p(x_0,\xi_0)=0$. In this situation we have
$\RR\times \{(x_0,\xi_0)\} \subset \mathrm{Fix}_{E_c}$.

When $E$ is not critical and the periodic orbits satisfy a
condition of non-degeneracy, the asymptotics behavior of
Eq.(\ref{Def trace}) is well determined by the closed orbits of
$\Phi_t$ on $\Sigma_{E}$ and the geometry of $\Sigma_E$. For the
full treatment of this problem, and a complete formulation
of the asymptotic expansion, we refer to \cite{BU,PU}.\medskip\\
\textbf{Removing divergences.} We explain now shortly why the
problem stated in Eq.(\ref{Def trace}) leads to a mathematically
rigorous version of the Gutzwiller formula. First, for each $h>0$
the sum is finite and a fortiori convergent. A convenient choice
of $\varphi$ also ensures that this quantity has an asymptotic
expansion when $h\rightarrow 0$ independently from the choice of
$\varepsilon>0$ up to corrections of order $\mathcal{O}(h^\infty)$
as long as $\varepsilon$ stays strictly positive. Such a
difference of size $\mathcal{O}(h^\infty)$ plays no r\^ole since
the discussion will be based on some finite order asymptotics
w.r.t. $h$.

On the other side, only the periods of $\Phi_t$ inside
$\rm{supp}(\hat{\varphi})$, the support of the Fourier transform:
\begin{equation*}
\hat{\varphi}(t)=\int\limits_{\mathbb{R}} e^{itx}\varphi(x)dx,
\end{equation*}
contribute in the asymptotic expansion. This principle is useful
since when $\rm{supp}(\hat{\varphi})$ is compact then finitely
many closed orbits of $\Sigma_E$ contribute and the second
divergence is solved. Hence if $\hat{\varphi}\in
C_0^{\infty}(\mathbb{R})$, the space of smooth functions with
compact support, $\varphi$ is in the Schwartz space
$\mathcal{S}(\mathbb{R})$. Since elements of
$\mathcal{S}(\mathbb{R})$ are smooth with exponential decay at
infinity, no divergence occurs and the size of $\varepsilon$ is
irrelevant, up to an $\mathcal{O}(h^\infty)$-error.

Finally, in Eq.(\ref{Def trace}) the scaling w.r.t. $h$ is
important. With this choice and via Fourier transform
considerations, we can use the semi-classical propagator
$U_{h}(t)=\exp(itP_{h} /h)$, solution of the Schr\"odinger
equation:
\begin{equation*}
-ih \partial_t U_{h}(t)=P_h U_{h}(t),
\end{equation*}
to obtain a precise control w.r.t. $h$. Roughly, $U_{h}(t)$ can be
expanded w.r.t. $h$ via a so-called WKB approximation. This
expansion also provides the explicit relation with the classical
dynamics. The precise technical justifications are given in
section 3.
\subsection{Critical values and contributions of equilibria.}
In the previous section we heuristically outlined a relation valid when the semi-classical parameter tends to 0:
\begin{equation*}
\lim_{h\rightarrow 0}\Upsilon (E,h,\varphi )\rightleftharpoons
\mathrm{Fix}_E=\{(t,x,\xi)\in\mathbb{R}\times \Sigma _{E} \text{ /
} \Phi_{t}(x,\xi)=(x,\xi)\}.
\end{equation*}
Meaning that the asymptotic behavior of the left hand side can be
expressed in terms of distributions generated by fixed point of
the flow.

In the r.h.s any point $(x,\xi)$ of a periodic orbit $\gamma$
appears only at times $kT^\sharp_\gamma$, $k\in\mathbb{Z}^*$,
where $T^\sharp_\gamma$ is the primitive period of $\gamma$ orbit.
Also any point of the energy surface contributes for $t=0$ since
the flow is the identity at $t=0$. But an equilibrium
$(x_0,\xi_0)$ satisfies $\Phi_t(x_0,\xi_0)=(x_0,\xi_0)$ for all
$t$. Hence when $E$ is no more a regular value the nature of the
set of fixed point changes and some new contributions appear in
the asymptotic expansion. These new contributions can be qualified
of \textit{new wave invariants} (see below) and are extremely
important for the inverse spectral problem.

When $E=E_c$ is a critical value of the principal symbol $p$, the
asymptotic behavior of Eq.(\ref{Def trace}) is more complicated
and is closely related to the geometry of the flow inside
$\Sigma_{E_c}$. The presence of classical equilibria inside
$\Sigma_{E_c}$ and the stability of the flow near the critical set
affects strongly the nature of the asymptotic expansions. For a
non-degenerate critical point, i.e. when $d^{2}p(x_0,\xi_0)$ is an
invertible matrix when $dp(x_0,\xi_0)=0$, the reader can consult
\cite{BPU}. The problem is treated there for quite general
operators, also including the case of a manifold of critical
points, but for $\rm{supp}(\hat{\varphi})$ small around the
origin. For Schr\"{o}dinger operators on $\mathbb{R}^n$ and
$\mathrm{supp}(\hat{\varphi})$ compact but arbitrary, the results
of \cite{BPU} are improved in \cite{KhD1}.

Two important problems occur in presence of critical points.
First, at every point where $dp=0$ the surface $\Sigma_{E_c}$ and
the metric of $\Sigma_{E_c}$ are not smooth. Next, the
determination of the asymptotic expansion w.r.t. $h$ can be very
difficult. The point is that $\Upsilon(E,h,\varphi)$ can be
expressed in terms of oscillatory integrals:
\begin{equation*}
I(h)= \int\limits_{\mathbb{R}\times T^*\mathbb{R}^{n}} a(t,x,\xi)
e^{\frac{i}{h} \psi(t,x,\xi)}dtdxd\xi,\text{ } h \rightarrow 0^+.
\end{equation*}
The oscillating coefficient $h^{-1}$ is precisely imposed by the
scaling w.r.t. $h$ in Eq.(\ref{Def trace}) and plays an important
r\^ole since $I(h)$ oscillates fast in the semi-classical regime.
Via the WKB approximation, the phase $\psi$ is related to the flow
so that the asymptotic behavior of $I(h)$ is determined by the
closed orbits. The technical problem is that, in presence of an
equilibrium, $\psi$ has some degenerate critical points. The
stationary phase method cannot be applied and the asymptotic
expansion of $I(h)$ is radically different : e.g. some terms
$h^\alpha$, $\alpha\in\mathbb{Q}$ and powers of $\log(h)$
generally appear in this setting (see below). Also the nature of
these new terms can be very different since for example they can
be associated to some distributions acting on $\varphi$
with a continuous support (e.g. the full set of real numbers or a half-line).\medskip\\
\textbf{Wave-invariants.} A classical approach (used in section 3) is to study the asymptotic behavior, as $h\rightarrow 0^+$,
of the localized trace:
\begin{equation*}
\Omega(E,h,t)=\mathrm{Tr}\, \left(\Theta(P_h) e^{-\frac{it}{h}
(P_h-E)}\right), \Theta\in C_0^\infty.
\end{equation*}
I will follow now the terminology used in \cite{Hez}. Under
certain assumptions (see section 3 and 4), and for $E$ regular it
is well known that $\Omega$ admits an asymptotic expansion of the
form:
\begin{equation*}
\Omega(E,h,t)\sim \sum\limits_{j=-n}^{\infty} a_j(E,t) h^j, \text{ as } h\rightarrow 0^+.
\end{equation*}
The coefficients $a_j(E,t)$ are some tempered distributions on the
line $\mathbb{R}_t$ and are called \textit{wave invariants} of
$P_h$. These distributions have a different expression when $E$
varies and many of them are continuous functions of $E$ as long as
we do not cross critical levels of the energy function $p$.

When $E\rightarrow E_c$ ($E_c$ stands for critical levels) one can
observe a discontinuity in the asymptotic expansion. Also some new
coefficients generally appear, since we have the asymptotics:
\begin{equation}\label{formal expand}
\Omega(E_c,h,t)\sim \sum\limits_{k=0}^{n-1}\sum\limits_{j=-n_0}^{\infty} a_{j,k}(E_c,t) h^\frac{j}{p} \log(h)^k, \text{ as } h\rightarrow 0^+,
\end{equation}
for some $p\in \mathbb{N}^*$, see \cite{BPU,Cam0,Cam3,Cam4,KhD1}
for details and examples. These new coefficients are called
\textit{new wave invariants} and the top order coefficient of the
expansion in Eq.(\ref{formal expand}) contains many information on
the shape of the symbol.
\subsection{Results and strategy}
Our first objective is to relate some variations in the discrete
spectrum of $P_{h}$ with the presence of fixed points for the
classical system: this principle detects the presence of new wave invariants and a fortiori critical energy levels.
Secondly, we establish that the precise
knowledge of such a spectral fluctuation can describe the
singularity of the potential (or of the symbol for general operators).
In theory, such a determination is possible since the contributions of equilibriums are highly
sensitive to the local shape of $V$ and are extremely persistent when the test function $\varphi$ varies.

For Schr\"odinger operators, we will consider the case of a potential $V$ with finitely many
critical points $x_0^j$ attached to local homogeneous extremum of
$V$. An immediate consequence is that $p$ admits, locally, a
unique critical point $(x_0^j,0)$ on the surface
$\Sigma_{E_{c}^j}=\{(x,\xi)\in \mathbb{R}^{2n}\text{ / }
\xi^2+V(x)=V(x_0^j)\}$. A typical example is a polynomial double
well in dimension 1 where 3 critical points occur at the 2 minima
and at the maximum of V (see figure 1).
\begin{figure}[h!] \label{double well}
   \begin{minipage}[l]{.35\linewidth}
       \epsfig{figure=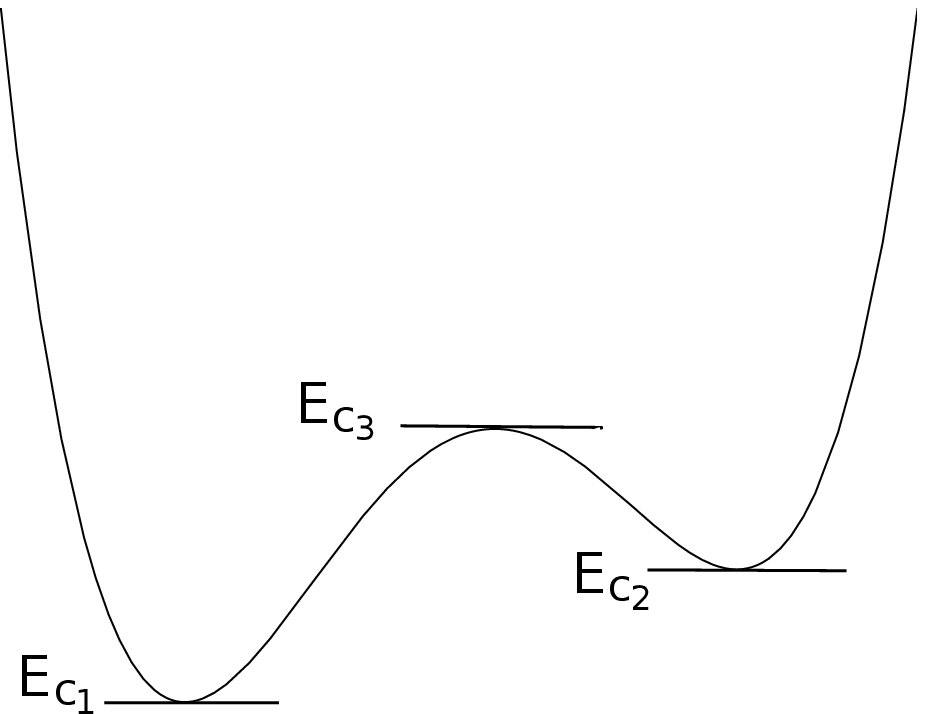,width=1. \textwidth}
   \end{minipage}
    \begin{minipage}[r]{.35\linewidth}
       \epsfig{figure=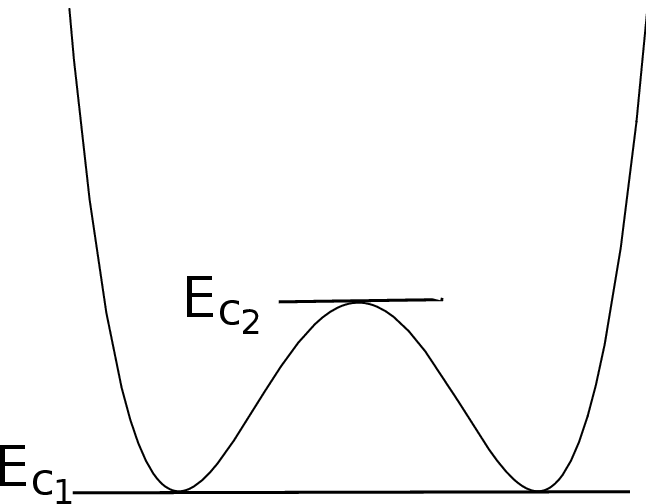,width=1. \textwidth}
   \end{minipage}
      \caption{Non-symmetric and symmetric double well.}
\end{figure}

In fact using certain generalizations of stationary phase methods,
necessary if the phase has some degenerate critical points one can
derive a very precise relation between the spectrum and the set of
fixed points of the flow also including the new wave invariants.
Once this relation is established in the form of asymptotic
expansion, the main results follow since:
\begin{itemize}
\item Equilibriums have a continuous contribution w.r.t. the time $t$.%
\item Shrinking $\mathrm{supp}(\hat{\varphi})$ erases all other
contributions: Weyl-terms and periodic orbits.%
\item The remaining contribution, given by new wave invariants,
displays some nice information on $V$.
\end{itemize}
The first assertion simply means that a fixed point
generally\footnote{Some discrete contributions can also sometimes
occur as pointed out in \cite{BPU} or \cite{Cam0}. But, to attain
our objectives, we can avoid to include them in the spectral
estimates.} contributes to the asymptotic expansion of $\Upsilon
(E^j_c,h,\varphi )$ in the form $h^\alpha \log(h)^\beta
\left\langle T_{\alpha,\beta},\hat{\varphi}\right\rangle$ where
$T_{\alpha,\beta}$ is a distribution such that
$\mathrm{supp}(T_{\alpha,\beta})=\mathbb{R}$, $[a,\infty]$ or
$[-\infty,a]$. Contrary to standard periodic orbits whom
contributions are supported in the set of periods, such a term
supported on the line cannot be erased just by shrinking the
support of $\hat{\varphi}$.

For example, if $\mathrm{supp}(\hat{\varphi})$ contains no period
of the flow our analysis follows if we view $\Upsilon(E,h,\varphi)$ as a function of $E$ :
\begin{itemize}
\item The order w.r.t $h$ of $\Upsilon(E,h,\varphi)$ changes when
$E\rightarrow E^j_c$ (Prop. \ref{Main}). \item This discontinuity
at $E^j_c$ describes the shape of $V$.
\end{itemize}
This indicates the presence of an equilibrium for $\Phi_t$, a
fortiori of a critical point for $V$.
\begin{remark}
\rm{For a degenerate singularity the information is more difficult
to interpret compared to a non-degenerate singularity (see section
2).}
\end{remark}
\section{Hypotheses and main result.}
Let $p(x,\xi)=\xi^2 +V(x)$, where the potential $V$ is real valued
and smooth on $\mathbb{R}^n$. To this Hamiltonian is attached the
operator $P_{h}=-h^2\Delta+V(x)$ and by a classical result, see
\cite{Ber-Shu}, $P_{h}$ is essentially self-adjoint starting from
a dense domain of $L^2(\mathbb{R}^n)$ when $V$ is bounded from
below and with tempered growth.
\begin{remark}
\rm{In this section, we are here mainly interested in the case of
Schr\"odinger operators but generalizations to an $h$-admissible
operator (e.g. in the sense of \cite{[Rob]}) are given in section
5.}
\end{remark}
First, to obtain a well defined spectral problem, we use:\medskip\\
$(\mathcal{H}_{1})$ $V\in C^{\infty}(\mathbb{R}^n)$. \textit{There exists } $C\in \mathbb{R}$ \textit{ such that }%
$\liminf\limits_{\infty}V >C$.\medskip\\
Note that $(\mathcal{H}_1)$ is always satisfied if $V$ goes to
infinity at infinity. Now, consider an energy interval
$J=[E_1,E_2]$ with $E_2<\liminf\limits_{\infty}V$. In the
following we note :
\begin{equation}
J(\varepsilon)=[E_1-\varepsilon,E_2+\varepsilon].
\end{equation}
For $\varepsilon<\varepsilon_0$ the set $p^{-1}(J(\varepsilon))$
is compact. By Theorem 3.13 of \cite{[Rob]} the spectrum $\sigma
(P_{h})\cap J(\varepsilon)$ is discrete and consists in a
sequence:
\begin{equation*}
\lambda _{1}(h)\leq \lambda _{2}(h)\leq ...\leq \lambda _{j}(h),
\end{equation*}
of eigenvalues of finite multiplicities, if $\varepsilon$ and $h$
are small enough. In general such a condition that the pullback of
$J$ or $J(\varepsilon)$ by $p$ is compact is sufficient to obtain
a discrete spectrum. This is not necessary as shows the
non-confining potential $x^2y^2$ on $\mathbb{R}^2$, see
\cite{Sim2} where precise spectral estimates are given for such
potentials.

The central object of study is the spectral distribution:
\begin{equation}
\Upsilon (E,h,\varphi)=\sum\limits_{\lambda _{j}(h)\in
J(\varepsilon)}\varphi (\frac{\lambda _{j}(h)-E}{h}),
\label{Objet trace}
\end{equation}
and, more precisely, the asymptotic information contained in this
object as $h\rightarrow 0^+$. To avoid any problem of convergence
we impose the condition: \medskip\\
$(\mathcal{H}_{2})$\textit{ We have }$\hat{\varphi}\in C_{0}^{\infty
}(\mathbb{R})$ \textit{ with a sufficiently small support near the
origin.}
\begin{remark}
\rm{$(\mathcal{H}_{2})$ is used to erase contributions of non-trivial closed orbits and can be
relaxed to $\hat{\varphi}\in C_{0}^{\infty }(\mathbb{R})$ with a
weaker result. A more precise description of
$\mathrm{supp}(\hat{\varphi})$ is given in Lemma \ref{periods}.
For a non-degenerate minimum, it is more comfortable to assume
that $\mathrm{supp}(\hat{\varphi})$ contains no period of
$d\Phi_t(z_0)$. Some singularities, not directly relevant here, are
generated by these periods and we refer to \cite{BPU,KhD1} for a
detailed study of these contributions.}
\end{remark}
To simplify notations we write $z=(x,\xi)\in \mathbb{R}^{2n}$ and
$\Sigma_{E}=p^{-1}(\{E\})$  and we use the subscript $E_c$ to
distinguish out critical values of $p$. Of course one can also
work with
$T^*\mathbb{R}^n\simeq\mathbb{R}^{n}\times\mathbb{R}^{n}$. In $J$
there is finitely many critical values $E_c^1,...,E_c^l$ and in
$p^{-1}(J)$ finitely
many fixed points $z_0^1,...,z_0^m$, $m\geq l$. We impose now the type of singularity:\medskip\\
$(\mathcal{H}_{3})$\textit{ On each }$\Sigma _{E_c^j}$ \textit{the
symbol }$p$\textit{ has isolated critical points
}$z_{0}^j=(x_{0}^j,0).$ \textit{These critical points can be
degenerate but are associated to a local extremum of $V$:
\begin{equation}\label{form pot}
V(x)=E_c+ V_{2k}(x)+\mathcal{O}(||x-x^j_0||^{2k+1}), \text{
}k\in\mathbb{N}^{*},
\end{equation}
where $V_{2k}$, homogeneous of degree $2k$, is definite positive
or negative.}

\begin{remark} \label{critical point}
\rm{For non-degenerate singularities we can apply the results of
\cite{BPU,Cam0,KhD1} and the extremum condition is not really
necessary. We will recall and use these results in the next
section. But for a degenerate critical point of $V$ the extremum
condition is required since, to our knowledge, the contribution of
such a singularity to the wave expansion is unknown.}
\end{remark}
The next assumption, erases the mean values, i.e. the heat-like
invariants or so-called Weyl-terms, in the trace formula:\medskip\\
$(\mathcal{H}_4)$ \textit{$\hat{\varphi}$ is flat at 0, i.e. $\hat{\varphi}^{(j)}(0)=0$, $\forall j\in \mathbb{N}$}.\medskip\\
We can weaken condition $(\mathcal{H}_4)$ to
$\hat{\varphi}^{(j)}(0)=0$, $\forall j\leq k_0$, where $k_0\in
\mathbb{N}^*$ depends only on the degree of the singularities of
$V$ (see section 4) without essential change. Such a function
$\varphi$ exists and is easy to construct. Pick $\phi\in
C_0^{\infty}(\mathbb{R})$, $\mathrm{supp}(\phi)\subset [-M,M]$,
then $\hat{\varphi}(t)=t^{2j_0}\phi(t)$ satisfies our hypotheses.
In this case, we can chose the function $\phi$ even so that
$\varphi$ is real.

Finally, to relax a bit $(\mathcal{H}_2)$ we need a control on
the contribution of closed orbits. To do so, we impose the classical condition :\medskip\\
$(\mathcal{H}_5)$ \textit{ All periodic trajectories of the flow are
non-degenerate.}\medskip\\
Non-degenerate closed orbits are those whose Poincar\'e map does
not admit 1 as eigenvalue and are isolated. The non-degeneracy
condition on orbits is not a central argument in this work and is
only used to control the order w.r.t. $h$ of the contribution of
closed orbits (i.e. the order of the usual wave-invariant). One
could also impose a condition of 'clean-flow' to consider
families/submanifolds or bunches of closed orbits of positive
dimension. These conditions on the classical dynamics can be
simply discarded as soon as we have a strictly positive
lower-bound on periods of closed-orbits like in Lemma \ref{periods}.\medskip\\
\textbf{Detecting critical levels.} The first result shows how to
detect critical energy levels by revealing a singularity in the
spectral estimates:
\begin{proposition}[\bf{Spectral variation}]\label{Main}$\,$\\
Assume that conditions $(\mathcal{H}_{1})$ to $(\mathcal{H}_{4})$
are satisfied. As $h$ tends to $0^+$, we have:
\begin{equation*}
\Upsilon(E,h,\varphi)= \left\{
\begin{matrix}
\mathcal{O}(h^\infty) \text{ if } E\in [E_1,E_2]\backslash \{E_{c}^1,..., E_{c}^l\},\\
\mathcal{O}(f_j(h)) \text{ if } E=E_{c}^j,\text { }j\in
\{1,...,l\},
\end{matrix}\right.
\end{equation*}
where each $f_j(h)$ has a finite order w.r.t. $h$.
\end{proposition}
For a non-degenerate critical point the coefficients $f_j$ can be
determined explicitly. It is in general possible to predict a full
asymptotic expansion but to obtain an invariant formulation of all
the coefficients can be difficult since for degenerate critical
points the method you have to use is more complicated than the
usual stationary-phase formula.

For example, if $E=E_{c}^j$ and the surface $\Sigma_{E_c^j}$
carries a single minimum of degree $2k$, $k>1$, using the main
result of \cite{Cam3} we obtain:
\begin{equation} \label{report1}
f_j(h)=D(n,k,\varphi) h^{\frac{n}{2}+\frac{n}{2k}-n}.
\end{equation}
In the same situation, but for a local maximum of $V$, using the
results of \cite{Cam4}, we can obtain a logarithm of $h$:
\begin{equation}\label{report2}
f_j(h)=D(n,k,\varphi)
h^{\frac{n}{2}+\frac{n}{2k}-n}\log(h)^{j},\text{ } j=0 \text{ or }
1.
\end{equation}
In Eq.(\ref{report1}) and Eq.(\ref{report2}) the coefficient $D$
is a tempered distributions acting $\varphi$ characteristic from
the nature of the critical point. In fact if the critical surface
carries more than one critical point then $f_j$ is the sum of
their respective contributions. Note that for $n=1$ and $k>1$ the
singular term has negative order w.r.t. $h$. A detailed
formulation of the coefficients $f_j(h)$ is given in Propositions
\ref{minimum},\ref{maximum1} and
\ref{maximum2}.\medskip\\
\textbf{Inverse spectral results.} An interesting property is that
in the singularity of $\Upsilon(E,h,\varphi)$ when $E\rightarrow
E_c$ the order w.r.t. $h$ but also the constants of the top
order-coefficients, describes partially the shape of $V$:
\begin{theorem}[\bf{Inverse result for Morse-critical points}]\label{Main1}$\,$\\
Assume that $\Sigma_{E_c}$ carries exactly one critical point associated to a non-degenerate critical point $x_0$ of $V$. Then the discontinuity of
$\Upsilon(s,h,\varphi)$ at $s=E_c$ determines the spectrum of $d^2V(x_0)$.
\end{theorem}
This result follows from the special form of the
Duistermaat-Guillemin-Uribe density at a non-degenerate critical
point. Observe that in particular we retrieve the Morse index of
$V$ at $x_0$ (number of positive eigenvalues minus number of
negative eigenvalues of $d^2V(x_0)$). Once more, because of the
invariance under coordinates permutations, or under a rotation of
the potential around the critical point $x_0$, it is in general
not
possible to retrieve the quadratic form $d^2 V(x_0)$ in a given system of coordinates $x$.\medskip\\
For a degenerate homogeneous singularity we have also a nice result:
\begin{theorem}[\bf{Inverse result. Degenerate critical points of $V$.}]\label{Main2}$\,$\\
Assume that there is exactly one critical point $(x_0,0)$ on the
singular energy surface $\Sigma_{E_c}$. Assume that $x_0$ is
attached to an homogeneous maximum or minimum of the potential of
degree $2k$. Then the discontinuity of $\Upsilon(E,h,\varphi)$ at
$E=E_c$ determines:
\begin{itemize}
\item The degree $2k$ of the critical point of $V$.%
\item The spherical mean-value of the germ of $V$ in $x_0$:
\begin{equation*}
A(V)=\int\limits_{\mathbb{S}^{n-1}} |V_{2k}(\theta)|^{-\frac{n}{2k}} d\theta.
\end{equation*}
\end{itemize}
\end{theorem}
Observe that $A(V)$ itself is invariant under rotation and
translation of $V$. Both results of Theorems \ref{Main1} and
\ref{Main2} are limited in presence of multiple equilibriums on
the same surface since the sum of contributions of each critical
point could lead to a compensation or to several contributions of
exactly same nature and order. In general if $V$ is not a Morse
function, or if a surface of energy $\Sigma_{E_c}$ carries more
than one critical point, eigenfunction-estimates seems to be
required to have a well-posed inverse problem. Finally, at the end
of the article we will show up some interesting invariants for
pseudo-differential operators with homogeneous singularities
attached to extremum of the symbol.
\begin{remark} \rm{I would like to emphasize that
a maximum is more difficult to detect contrary to a local minimum
which is an isolated point of the energy surface (locally
$\Sigma_{E_c}$ is just a point). A similar result holds for an
operator of the form $T(\xi)+V(x)$ (kinetic plus potential energy)
where $T$ is convex near the origin and $V$ has a local minimum at
$x_0$. Moreover a maximum, attached to an unstable critical point
of the flow, is much more complicated to treat with semi-classical
methods. See \cite{Cam4} for a detailed study.}
\end{remark}
\noindent\textbf{Longer range estimates.} In $(\mathcal{H}_2)$ the
condition that $\mathrm{supp}(\hat{\varphi})$ is small implies a
very accurate spectral estimate (e.g. by a Paley-Wiener estimates
for the decay of $\varphi$). It is possible to relax this
assumption but the result is a bit weaker:
\begin{corollary}\label{weak}
Assume that conditions $(\mathcal{H}_1)$, $(\mathcal{H}_3)$,
$(\mathcal{H}_4)$ and $(\mathcal{H}_5)$ are satisfied and that
$\hat{\varphi}\in C_0^{\infty}(\mathbb{R})$, then we obtain:
\begin{equation*}
\Upsilon(E,h,\varphi)= \mathcal{O}(1) \text{ if } E\in
[E_1,E_2]\backslash \{E_c^1,..., E_c^l \}.
\end{equation*}
For critical values of $p$, estimates are the same as in
Proposition \ref{Main}.
\end{corollary}
\noindent The justification, see section 4, is that in this case
the asymptotics is given by a finite sum over periodic orbits of
energy $s$. This result is weak if the singularity of $V$ is
non-degenerate since the equilibrium has a contribution of degree
0 w.r.t. $h$, see Propositions
\ref{minimum},\ref{maximum1},\ref{maximum2} or section 3 of
\cite{BPU}.

In theory there is always a variation when $E\rightarrow E_c$ but
of course this effect can be harder to detect if there is no
change in the order w.r.t. $h$. In that situation there is only a
discontinuity in the top-order coefficient w.r.t. $h$ so that the
result can be qualified of 'weaker'.
\section{Oscillatory representation.}
The construction below is more or less classical and will be
sketchy. The only change with the usual construction, around a
single energy level, is that we use a more global localization
around $J=[E_1,E_2]$. Strictly speaking, with $(\mathcal{H}_1)$,
we could also consider $]-\infty, E_2]$ since there is no
eigenvalue below a fixed energy level $E_0$ given by the minimum
of the quadratic form attached to $P_h$. Let be $\varphi \in
\mathcal{S}(\mathbb{R})$ with $\hat{\varphi}\in C_{0}^{\infty
}(\mathbb{R})$, we recall that:
\begin{equation*} \Upsilon (E,h,\varphi)
=\sum\limits_{\lambda _{j}(h)\in J(\varepsilon)}\varphi
(\frac{\lambda _{j}(h)-E}{h}), \text { }
J(\varepsilon)=[E_1-\varepsilon,E_2+\varepsilon],
\end{equation*}
with $p^{-1}(J(\varepsilon))$ compact in $T^{\ast
}\mathbb{R}^{n}$. For $\varepsilon>0$ small enough, we localize
around $J$ with a cut-off $\Theta \in C_{0}^{\infty
}(]E_1-\varepsilon ,E_2+\varepsilon \lbrack )$, such that $\Theta
=1$ on $J$ and $0\leq \Theta \leq 1$ on $\mathbb{R}$. We
accordingly split-up our spectral distribution as:
\begin{equation*}
\Upsilon (E,h,\varphi)=\Upsilon _{1}(E,h,\varphi)+\Upsilon
_{2}(E,h,\varphi),
\end{equation*}
with :
\begin{gather*}
\Upsilon _{1}(E,h,\varphi)=\sum\limits_{\lambda _{j}(h)\in
J(\varepsilon)}(1-\Theta )(\lambda _{j}(h))\varphi
(\frac{\lambda _{j}(h)-E}{h}),\\
\Upsilon _{2}(E,h,\varphi)=\sum\limits_{\lambda _{j}(h)\in
J(\varepsilon)}\Theta (\lambda _{j}(h))\varphi (\frac{\lambda
_{j}(h)-E}{h}).
\end{gather*}
Since $\varphi\in \mathcal{S}(\mathbb{R})$ a classical estimate,
see e.g. Lemma 1 of \cite{Cam1}, is:
\begin{equation}
 \Upsilon _{1}(E,h,\varphi)=\mathcal{O}(h^{\infty }),
\text{ as } h\rightarrow 0^{+}.\label{S1(h)=Tr}
\end{equation}
By inversion of the Fourier transform we have:
\begin{equation*}
\Theta (P_{h})\varphi (\frac{P_{h}-E}{h})=\frac{1}{2\pi}\int\limits_{%
\mathbb{R}}e^{i\frac{tE}{h}}\hat{\varphi}(t)\mathrm{exp}(-\frac{it}{h}%
P_{h})\Theta (P_{h})dt.
\end{equation*}
The trace of the left hand-side is $\Upsilon _{2}(E,h,\varphi)$
and Eq.(\ref{S1(h)=Tr}) provides :
\begin{equation}\label{Trace S2(h)}
\Upsilon (E,h,\varphi)=\frac{1}{2\pi }\mathrm{Tr}\int\limits_{\mathbb{R}}e^{i%
\frac{tE}{h}}\hat{\varphi}(t)\mathrm{exp}(-\frac{it}{h}P_{h})\Theta
(P_{h})dt+\mathcal{O}(h^{\infty }).
\end{equation}
Eq.(\ref{Trace S2(h)}) is very close to the classical Poisson
summation formula on $\mathbb{S}^1$, see \cite{Sj-Zw} for a
discussion and an interpretation, since the r.h.s. is expressed
below in term of the classical dynamics. This asymptotic relation
justifies the terminology of \textit{trace formula}.

Moreover, the formulation in Eq.(\ref{Trace S2(h)}) shows that the
scaling w.r.t. $h$, imposed in the initial definition of
$\Upsilon(E,h,\varphi)$, is the best one since we will solve now
the semi-classical propagator homogeneously w.r.t. $h$. Let
$U_{h}(t)=\mathrm{exp}(-\frac{it}{h}P_{h})$ be the quantum
propagator. We approximate $U_{h}(t)\Theta (P_{h})$ by a Fourier
integral operator (FIO) depending on $h$. Let $\Lambda$ be the
Lagrangian manifold associated to the flow of $p$:
\begin{equation*}
\Lambda =\{(t,\tau ,x,\xi ,y,\eta )\in T^{\ast }\mathbb{R}\times T^{\ast }%
\mathbb{R}^{n}\times T^{\ast }\mathbb{R}^{n}:\tau =p(x,\xi
),\text{ }(x,\xi )=\Phi _{t}(y,\eta )\},
\end{equation*}
and $I(\mathbb{R}^{2n+1},\Lambda )$ the class of oscillatory
integrals based on $\mathbb{R}^{2n+1}$ and whose Lagrangian
manifold is $\Lambda$. The next result is a semi-classical version
of a well known result on the propagator, see e.g. Duistermaat
\cite{DUI1}.
\begin{theorem}
The operator $U_{h}(t)\Theta (P_{h})$ is an $h$-FIO associated to
$\Lambda$. For each $N\in\mathbb{N}$ there exists $U_{\Theta
,h}^{(N)}(t)$ with integral kernel in H\"ormander's class
$I(\mathbb{R}^{2n+1},\Lambda )$ and $R_{h}^{(N)}(t)$ bounded, with
a $L^{2}$-norm uniformly bounded for $0<h\leq 1$ and $t$ in a
compact subset of $\mathbb{R}$, such that:
\begin{equation*}
U_{h}(t)\Theta(P_{h})=U_{\Theta,h}^{(N)}(t)+h^{N}R_{h}^{(N)}(t).
\end{equation*}
\end{theorem}
This result provides the existence of an asymptotic expansion in
power of $h$ with a remainder that can be controlled since
$\mathrm{supp}(\hat{\varphi})$ is a compact. After perhaps a
reduction of $\varepsilon$, this remainder $R_{h}^{(N)}(t)$ is
estimated via:
\begin{corollary}
Let $\Theta _{1}\in C_{0}^{\infty }(\mathbb{R})$, with $\Theta
_{1}=1$ on $\mathrm{supp}(\Theta )$ and $\mathrm{supp}(\Theta
_{1})\subset ]E_1-2\varepsilon, E_2+2\varepsilon[$, then $\forall
N\in \mathbb{N}$:
\begin{equation*}
\mathrm{Tr}(\Theta (P_{h})\varphi (\frac{P_{h}-E}{h}))=\frac{1}{2\pi }%
\mathrm{Tr}\int\limits_{\mathbb{R}}\hat{\varphi}(t)e^{\frac{i}{h}%
tE}U_{\Theta ,h}^{(N)}(t)\Theta
_{1}(P_{h})dt+\mathcal{O}(h^{N-n}).
\end{equation*}
\end{corollary}
For a proof of this result, based on the cyclicity of the trace
and a priori estimates on the spectral projectors (see
\cite{[Rob]}), we refer to \cite{Cam1}. For the particular case of
a Schr\"odinger operator the BKW ansatz shows that the integral
kernel of $U_{\Theta,h}^{(N)}(t)$ can be recursively constructed
as:
\begin{gather*}
K_{h}^{(N)}(t,x,y)=\frac{1}{(2\pi h)^n}
\int\limits_{\mathbb{R}^n} b_{h}^{(N)}(t,x,y,\xi)
e^{\frac{i}{h} (S(t,x,\xi)-\left\langle y,\xi
\right\rangle)} d\xi,\\
b_{h}^{(N)}=b_0+h b_1+...+h^N b_N,
\end{gather*}
where $S$ satisfies the Hamilton-Jacobi equation:
\begin{equation*}
p(x,\partial _{x}S(t,x,\xi))+\partial _{t}S(t,x,\xi )=0,
\end{equation*}
with initial condition $S(0,x,\xi)=\left\langle x,\xi
\right\rangle$. In particular we obtain that:
\begin{equation*}
\{(t,\partial _{t}S(t,x,\eta ),x,\partial _{x}S(t,x,\eta
),\partial _{\eta }S(t,x,\eta ),-\eta )\}\subset \Lambda ,
\end{equation*}
and that the function $S$ is a generating function of the flow,
i.e.:
\begin{equation}
\Phi _{t}(\partial _{\eta }S(t,x,\eta ),\eta ) =(x,\partial
_{x}S(t,x,\eta )). \label{Gene}
\end{equation}
We insert this approximation in Eq.(\ref{Trace S2(h)}), we set
$x=y$ and we integrate w.r.t. $x$. Modulo an error
$\mathcal{O}(h^{N-n})$, we obtain that $\Upsilon(E,h,\varphi)$
equals:
\begin{equation}
\frac{1}{(2\pi h)^n}\int\limits_{\mathbb{R}\times
T^*\mathbb{R}^n}e^{\frac{i}{h}(S(t,x,\xi )-\left\langle x,\xi
\right\rangle +tE)}a_{h}^{(N)} (t,x,\xi
)\hat{\varphi}(t)dtdxd\xi, \label{gamma1 OIF}
\end{equation}
where $a_{h}^{(N)}(t,x,\eta )=b_{h}^{(N)}(t,x,x,\eta )$.
\begin{remark}
\rm{By Theorem 3.11 \& Remark 3.14 of \cite{[Rob]}, $\Theta
(P_{h})$ is $h$-admissible. Moreover, the symbol is compactly
supported in $p^{-1}([E_1-\varepsilon,E_2+\varepsilon])$. This
global result w.r.t. $E\in [E_1,E_2]$ allows to consider below
only oscillatory integrals with compact support for the evaluation
of the spectral distributions.}
\end{remark}
\noindent\textbf{Micro-localization of the trace.}\\
If $\psi\in C_0^{\infty}(T^{*}\mathbb{R}^{n})$, we recall that
$\psi^w(x,h D_x)$ is the linear operator obtained by
Weyl-quantization of $\psi$ and semi-classical quantization of
$\xi$. This means:
\begin{equation*}
\psi ^{w}(x,h D_{x})f(x)=\frac{1}{(2\pi h)^n}
\int\limits_{\mathbb{R}^{2n}} e^{\frac{i}{h} \left\langle
x-y,\xi\right\rangle} \psi(\frac{x+y}{2},\xi) f(y)dyd\xi.
\end{equation*}
Mainly, the contribution of an equilibrium $z_0\in \Sigma_{E_c}$
can be reached via:
\begin{equation}\label{trace local}
\Upsilon_{z_0}(E_c,h ,\varphi)=\frac{1}{2\pi }\mathrm{Tr}\int\limits_{\mathbb{R}}e^{i%
\frac{tE_{c}}{h}}\hat{\varphi}(t)\psi ^{w}(x,h D_{x})\mathrm{exp}(-\frac{i}{h}%
tP_{h})\Theta (P_{h})dt,
\end{equation}
where $\psi\in C_0^{\infty} (T^* \mathbb{R}^n)$ is equal to 1 near
$z_0$. Since the trace is a cyclic operation, we can use the
symbolic-calculus and to insert an $L^2$-bounded observable (here
a cut-off in the phase space) is in general not very expensive in
the FIO construction. This construction with smooth cut-off allows
to work with pseudo-differential partition of unity. This approach
is useful to obtain a weak generalization of our results in
presence of multiple equilibriums.

We recall some basic results on the symbolic calculus with FIO.
H\"ormander's class of distributions with Lagrangian manifold
$\Lambda$ over $\mathbb{R}^n$ is noted $I(\mathbb{R}^{n},\Lambda
)$. If $(x_{0},\xi _{0})\in \Lambda $ and $\varphi (x,\theta )\in
C^{\infty }(\mathbb{R}^{n}\times \mathbb{R}^{N})$ parameterizes
$\Lambda $ in a sufficiently small neighborhood $U$ of $(x_{0},\xi
_{0})$, then for each $u_{h}\in I(\mathbb{R}^{n},\Lambda )$ and
$\chi \in C_{0}^{\infty }(T^{\ast }\mathbb{R}^{n})$,
$\rm{supp}(\chi )\subset U,$ there exists a
sequence of amplitudes $c_{j}(x,\theta )\in C_{0}^{\infty }(\mathbb{R}%
^{n}\times \mathbb{R}^{N})$ such that for all $L\in\mathbb{N}$:
\begin{equation*}
\chi ^{w}(x,h D_{x})u_{h}=\sum\limits_{-d\leq j<L}h^{j}I(c_{j}e^{\frac{i}{h}%
\varphi })+\mathcal{O}(h^{L}).
\end{equation*}
Hence, for each $N\in \mathbb{N}^{*}$ and modulo an error
$\mathcal{O}(h^{N-d})$, the localized trace
$\Upsilon_{z_0}(E_c,h,\varphi)$ of Eq.(\ref{trace local}) can be
written as:
\begin{equation}
(2\pi h)^{-d}\int\limits_{\mathbb{%
R\times R}^{2n}}e^{\frac{i}{h}(S(t,x,\xi )-\left\langle x,\xi
\right\rangle +tE_c)}\tilde{a}_{h}^{(N)}(t,x,\xi
)\hat{\varphi}(t)dtdxd\xi . \label{gamma1 OIF2}
\end{equation}
To obtain the exact power $-d$ of $h$, we apply results of
Duistermaat \cite{DUI1} on the order of FIO. Since:
\begin{itemize}
\item $h$-pseudo-differential operators $\psi^w(x,h D_x)$ are of
order 0 w.r.t. $1/h$.%
\item the order of $U_{h}(t)\Theta (P_{h})$ is $-\frac{1}{4}$.
\end{itemize}
if we identify the operator with it's distributional-kernel we
have:
\begin{equation*}
\psi ^{w}(x,h D_{x})U_{h}(t)\Theta (P_{h})\sim (2\pi
h)^{-n}\int\limits_{\mathbb{R}^{n}}\tilde{a}_{h}^{(N)}(t,x,y,\eta )e^{\frac{i}{h}%
(S(t,x,\eta )-\left\langle y,\eta \right\rangle )}d\eta.
\end{equation*}
Multiplying by $\hat{\varphi}(t)e^{\frac{i}{h}tE_{c}}$ and passing
to the trace we find Eq.(\ref{gamma1 OIF2}) with $d=n$ and we
write again $\tilde{a}_{h}^{(N)}(t,x,\eta )$ for the diagonal
evaluation $\tilde{a}_{h}^{(N)}(t,x,x,\eta )$. In particular:
\begin{equation}
\tilde{a}_{h}^{(0)}(t,x,x,\eta )= \psi(x,\eta)
a_{0}(t,x,x,\eta),
\end{equation}
is independent of $h$ and is compactly supported w.r.t. $(x,\eta)$
since $\psi\in C_0^\infty(T^*\mathbb{R}^n)$.
\section{Proof of the main result.}
Let be $E_c$ any critical value in $[E_1,E_2]$ and $z_0$ an
equilibrium of $\Sigma_{E_c}$. We choose a function $\psi \in
C_{0}^{\infty }(T^{\ast }\mathbb{R}^{n})$, with $\psi =1\text{
near }z_{0}$, hence:
\begin{gather*}
\Upsilon _{2}(E_{c},h,\varphi) =\frac{1}{2\pi }\mathrm{Tr}\int\limits_{\mathbb{R}}e^{i%
\frac{tE_{c}}{h}}\hat{\varphi}(t)\psi ^{w}(x,h D_{x})\mathrm{exp}(-\frac{i}{h}%
tP_{h})\Theta (P_{h})dt\\
+\frac{1}{2\pi }\mathrm{Tr}\int\limits_{\mathbb{R}}e^{i\frac{tE_{c}}{h}}\hat{%
\varphi}(t)(1-\psi ^{w}(x,h
D_{x}))\mathrm{exp}(-\frac{i}{h}tP_{h})\Theta (P_{h})dt.
\end{gather*}
If there is no other singularity on $\Sigma_{E_c}$ with
$(\mathcal{H}_5)$ the asymptotic expansion of the second term is
given by the semi-classical trace formula on a regular level. For
finitely many critical point on $\Sigma_{E_c}$, we can repeat the
procedure. The first term is micro-local and precisely generates
the singularity in Theorem \ref{Main}. We note $\Omega$ the
discrete set of
critical points $z_0^j$ in $p^{-1}(J)$. \medskip\\
\textbf{Classical dynamics near the equilibrium.}\\
A generic critical points of the phase function of Eq.(\ref{gamma1
OIF}) satisfies the relations:
\begin{equation*}
\left\{
\begin{array}{c}
E=-\partial _{t}S(t,x,\xi ), \\
x=\partial _{\xi }S(t,x,\xi ), \\
\xi =\partial _{x}S(t,x,\xi ),
\end{array}
\right. \Leftrightarrow \left\{
\begin{array}{c}
p(x,\xi )=E, \\
\Phi _{t}(x,\xi )=(x,\xi ).
\end{array}
\right.
\end{equation*}
The right hand side coincide with the set $\mathrm{Fix}_E$ of
fixed points of $\Phi_t$ defined in section 1. This set generally
consist of closed trajectories of the flow inside $\Sigma_{E}$,
the energy surface (for $t=0$ the flow is the identity) and
finally equilibria/critical-points. By the non-stationary phase
lemma, outside of the critical set $\mathrm{Fix}_E$ the
contribution is of order $\mathcal{O}(h^{\infty})$.\medskip

To apply the stationary-phase methods it is important to study the
nature of the phase function along the critical-set
$\mathrm{Fix}_E$. For the regularity of the Hessian of the phase
function, the next lemma is particulary useful. We will often
denote points $(x,\xi)\in T^*\mathbb{R}^n$ of the phase space by a
single letter $z$.
\begin{lemma}\label{Lien phase/flot}
Let us define $\Psi (t,x,\xi )=S(t,x,\xi
)-\left\langle x,\xi \right\rangle +tE_{c},$ then if $z_{0}$ is
critical point of $\Psi ,$ we have the equivalence:
\begin{equation*}
d_{z}^{2}\Psi (z_{0})\delta z=0\Leftrightarrow d_{z}\Phi _{t}(z_{0})\delta
z=\delta z,\, \forall \delta z\in T_{z_{0}}(T^{\ast}\mathbb{R}^{n}).
\end{equation*}
In other words, degenerate directions of the phase
correspond to fixed points of the linearized flow at $z_0$.
\end{lemma}
The proof is standard and can for example be found in
\cite{KhD1,Cam0}. We recall that the linearized flow $d\Phi_t$ is
the differential of the flow $\Phi_t$ w.r.t. initial conditions
$z=(x,\xi)$. If we use Lemma \ref{Lien phase/flot}, we obtain for
our phase function
\begin{corollary}
\label{theo periode flot linéarise} A critical point
$(T,x_{0},\xi_{0})$ of $\Psi (t,x,\xi )$ is degenerate with
respect to $(x,\xi )$ if and only if $T$ is a period of the
linearized flow $d_{x,\xi }\Phi _{t}(x_{0},\xi _{0})$.
\end{corollary}
The next result is also well known, see, e.g., \cite{A-M}, from classical mechanics and differential geometry:
\begin{lemma} \label{formule flot linéarisé}
If $\partial_{x}p(x_{0},\xi _{0})=\partial _{\xi }p(x_{0},\xi _{0})=0$ then
$d_{x,\xi }\Phi _{t}(x_{0},\xi _{0})$ is the Hamiltonian flow of
the quadratic form $\frac{1}{2}d^{2}p(x_{0},\xi _{0})$ on
$T_{x_{0},\xi _{0}}(T^{\ast }\mathbb{R}^{n})$.
\end{lemma}
Hence, when $z_0$ is a critical point of $p$, the linear map
$w\mapsto d\Phi_t(z_0)w$ can be interpreted as the Hamiltonian
flow of $u \mapsto Q(u)=\langle \frac{1}{2}d^2p(z_0)u,u\rangle$.
Observe that for a Schr\"odinger operator a critical point is
always of the form $z_0=(x_0,0)$ with $dV(x_0)=0$ and the
quadratic form is:
\begin{equation*}
Q(u)=Q(u_1,u_2)=||u_1||^2 +\frac{1}{2} \langle d^2 V(x_0)u_2,
u_2\rangle, \, u_1\in \mathbb{R}^n,\, u_2\in \mathbb{R}^n.
\end{equation*}
A fortiori the map $d\Phi_t(z_0)$ is governed
by a quadratic Schr\"odinger operator.\\
\textbf{Degenerate critical point.} If the kernel of $d^2 V(x_0)$
is not trivial and contains a non-zero vector $v$ the
linearized-flow in the direction of $v$ is the flow of the free
Laplace operator and we have $d\Phi_t(z_0)(v,0)=(v,0)$ for all
$t$. Hence it is never possible to apply the Morse-lemma, and a
fortiori the stationary phase method (and this for any time $t$).
But if we restrict the study to some homogeneous singularity for
$V$ it is possible to find a local diffeomorphism changing the
phase $\Psi$ into a polynomial function $R$ (i.e. a local
normal-form for $\Psi$). For these normal forms it is possible to
generalize the stationary-phase method:
\begin{equation*}
\int\limits_{\mathbb{R}^{2n+1}} e^{\frac{i}{h} R(w)} a(w)dw \sim
\sum\limits_{j,k} h^{\alpha_j} \log(h)^k c_{j,k}(a), \text{ as }
h\rightarrow 0^+.
\end{equation*}
We will not review all these geometric and analytic results
individually. We will refer to \cite{BPU,Cam1,Cam3,Cam4,Cam5} for
the reduction of the phase function and the asymptotic of the
resulting oscillatory integrals. I mention that it is complicated
to express all distributional coefficients $c_{j,k}(a)$
invariantly since $a$ depends on the local-diffeomorphism $\chi$
transforming $\Psi$ into $R$. The top-order coefficients are
relatively easy to express, the next one start to depend on the
derivatives of $\chi$. To express all $c_{j,k}(a)$ in terms of $V$
would improve the inverse spectral result to higher derivatives of
$V$ at the critical point.\medskip\\
\textbf{Non-degenerate critical point of $V$.} When $d^2V(x_0)$ is
invertible the have a classical result. It is easy to check that
the sign of the quadratic map $Q$ determines the 'stable' and
'unstable' directions after diagonalizing $d^2V$ via an orthogonal
linear transformation.

Combining the results of Lemma \ref{formule flot linéarisé} and
Corollary \ref{theo periode flot linéarise}, to achieve our goal
it will be sufficient to stay below the smallest positive period
of the linearized flow. Working in some suitable local coordinates
near $z_0$ we can assume that the quadratic form attached to the
potential is of the form:
\begin{equation*}
V(x)=E_c -\sum\limits_{j=1}^{r} \alpha_j x_j^2 +\sum\limits_{j=r+1}^{n} \alpha_j x_j^2+\mathcal{O}(||x||^3), \text{ } \alpha_j>0.
\end{equation*}
This can be achieved by a translation, a change of linear
coordinates (via an orthogonal matrix) and an eventual permutation
of coordinates. Observe that all these linear transformations, in
particular the action of the orthogonal matrix, leave the Laplace
operator invariant. Of course the variable attached to indices
$j=1,...,r$ will generate hyperbolic functions (and hence never
degenerate in the sense of Corollary \ref{theo periode flot
linéarise}).\medskip\\
\textbf{A lower bound on primitive periods.}\\
The next result provides a global information on the smallest positive primitive periods of the classical flow. This lemma will be used
to extract the new wave invariants appearing at a critical energy level.
\begin{lemma}\label{periods}
There exists $T>0$, depending only on $V$ and $J=[E_1,E_2]$, such
that $\Phi_t(z)\neq z$ for all $z\in p^{-1}(J)\backslash \Omega$
and all $t\in]-T,0[\cup ]0,T[$.
\end{lemma}
\textit{Proof.} If $H_p$ is our hamiltonian vector field and
$z=(x,\xi)$ we have:
\begin{equation*}
|| H_p(z_1)-H_p(z_2)||^2= 4||\xi_1-\xi_2||^2
+||\nabla_xV(x_1)-\nabla_xV(x_2)||^2.
\end{equation*}
When $z_1$ and $z_2$ are in the compact $p^{-1}(J)$ there exists
$b>0$ such that:
\begin{equation*}
||\nabla_xV(x_1)-\nabla_xV(x_2)||\leq b ||x_1-x_2||.
\end{equation*}
Hence, there exists $a>0$ such that:
\begin{equation*}
||H_p(z_1)-H_p(z_2)||\leq  a ||z_1-z_2||,\text { }\forall z_1,z_2
\in p^{-1}(J).
\end{equation*}
The main result of \cite{Yor} shows that any periodic orbit inside
$p^{-1}(J)$ has a period $\tau \geq 2\pi/a>0$. The lemma follows
with $T:=T(V,J)=2\pi/a$. $\hfill{\blacksquare}$
\begin{remark}
\rm{The result of \cite{Yor} is optimal (for the harmonic
oscillator the previous inequality becomes an equality). Note that
$T$ is decreasing if one increase the size of $J$. Lemma
\ref{periods} provides a total control on the r.h.s. of the trace
formula. If $\hat{\varphi}\in C_0^{\infty}(]0,T[)$, the only
contribution arises from the set $\{(t,z_0),\text{
}t\in\mathrm{supp}(\hat{\varphi})\}$, i.e. from the new wave
invariants.}
\end{remark}
Now, we restrict our attention to the singular contribution
generated by one critical point. We check now the condition of
non-degeneracy of our phase-function. As it was explained in
section 2, for a non-degenerate extremum a minor technical problem
could occur. If $x_0$ is a maximum of the potential $d\Phi_t(z_0)$
has no non-zero period which ends immediately the discussion. If
$x_0$ is a minimum $d\Phi_t(z_0)$ is elliptic with primitive
periods $(T_1,..,T_n)$ generated by the eigenvalues of $d^2
V(x_0)$. But the constant $b$ of Lemma \ref{periods} is certainly
bigger than the spectral radius of $d^2 V(x_0)$ and hence we have
the inequality $T<\min \{ T_1,..,T_n\}$. Following the approach of
\cite{BPU,Cam0,KhD1}, if $\mathrm{supp}(\hat{\varphi})\subset
]-T,T[$ the associated contribution is smooth on
$\mathrm{supp}(\hat{\varphi})\backslash \{0\}$.

For a degenerate critical point $z_0$ as in $(\mathcal{H}_3)$ a
surprising result, established in \cite{Cam3,Cam4}, is that the
only singularity, for the first new wave invariant, is located at
$t=0$. Hence the condition $\hat{\varphi}\in C_0^{\infty}(]0,T[)$
or $\hat{\varphi}\in C_0^{\infty}(]-T,0[)$ is sufficient to
determine the new-wave invariants generated by $z_0$.
\subsection{Non-degenerate critical points.}
With the previous considerations, we consider $z_0\in
\Sigma_{E_c}$ a non-degenerate critical point of $p$, a fortiori
isolated, and we assume that $\mathrm{supp}( \hat{\varphi})\subset
]-T,0[\cup ]0,T[$. In this setting we know that the phase function
$\Psi$ (introduced in Lemma \ref{Lien phase/flot}) has a non
degenerate critical point in $z_0$ for all $t\in \mathrm{supp}(
\hat{\varphi})$.

In our setting, up to a change of local coordinates, we can assume
that $p(x,\xi)=p_2(x,\xi)+\mathcal{O}(||x||^3$ (near the origin
$z_0=(0,0)$), where the quadratic form is:
\begin{equation}
p_{2}(x,\xi )=(\sum_{j=1}^{r}(\xi_j^2 -\alpha_j x_{j}^{2}))+(\sum_{j=r+1}^{n} (\xi _{j}^{2}+\alpha_j x_j^2)).
\label{truc}
\end{equation}
The flow of $p_{2}$, viewed as an element of $\mathrm{End}
(T_{0}(T^{\ast }\mathbb{R}^{n}))\simeq
\mathrm{End}(\mathbb{R}^{2n})$, is
\begin{gather*}
\mathrm{exp}(tH_{p_{2}})(x,\xi )=A(t)\left(
\begin{array}{l}
x \\
\xi
\end{array}
\right) ,\\%
A(t)=\left(
\begin{array}{cccc}
a(t) & 0 & e(t) & 0 \\
0 & b(t) & 0 & f(t) \\
c(t) & 0 & a(t) & 0 \\
0 & d(t) & 0 & b(t)
\end{array}
\right) ,
\end{gather*}
where $(x,\xi )=(x',x'',\xi',\xi'')$, $x',\xi' \in \mathbb{R}^r$,
$x'',\xi'' \in  \mathbb{R}^{n-r}$, and
\begin{equation*}
\left\{
\begin{array}{l}
a(t)=\mathrm{diag}(\,\mathrm{ch}(\alpha_{j}t)),\\
b(t)=\mathrm{diag}(\,\cos (\alpha_j t)),\\
c(t)=\mathrm{diag}(\, \alpha_j \mathrm{sh}(\alpha_j t)),\\
d(t)=\mathrm{diag}(\,-\alpha_j \sin (\alpha_j t)),\\
e(t)=\mathrm{diag}( \, \frac{1}{\alpha_j} \mathrm{sh}(\alpha_j t)),\\
f(t)=\mathrm{diag}(\, \frac{1}{\alpha_j} \sin (\alpha_j t)),
\end{array}
\right.
\end{equation*}
and the symbol "diag" means diagonal matrix. Observe that for $a(t),c(t),e(t)$ the index $j$ varies in $\{1,r\}$, for the others it varies in $\{r+1,n\}$
and that $\mathrm{det}{A}(t)=1$ for all $t$.\medskip\\
\textbf{Local reduction of the phase function of our FIO.} Since
$S$ is solution of the Hamilton-Jacobi equation we have (locally)
the relation:
\begin{equation*}
\Phi_t( \partial_\xi S(t,x,\xi),\xi)=(x,\partial_x S(t,x,\xi)).
\end{equation*}
But since $\Phi_t(z_0)=z_0$, writing a Taylor expansion for:
\begin{equation*}
S(t,x,\xi)=S_2(t,x,\xi)+\mathcal{O}(||(x,\xi)||^3,
\end{equation*}
and $\Phi_t$ (always near the critical point $z_0=0$) we get that $S_2$ satisfies the relation:
\begin{equation*}
d\Phi_t(0) (\partial_\xi S_2,\xi)+\mathcal{O}(||(x,\xi)||^2)= (x,\partial_x S_2)+\mathcal{O}(||(x,\xi)||^2),
\end{equation*}
where the linear map $d\Phi_t(0)$ can be determined as above. This linear system is regular exactly when $\det d(t)\neq 0$.
As a consequence, we pick:
\begin{equation*}
l=\inf (\frac{2\pi}{\alpha_j},\, j\in \{r+1,...,n\}),
\end{equation*}
and an interval $L\subset ]0,l[$. For each $t\in L$, it follows that we can find a (time-dependant) change of coordinates $(t,x,\xi)\mapsto (t,\chi)$,
well defined in a small neighborhood of $\mathrm{supp}( \hat{\varphi})\times \{z_0\}$ such that:
\begin{gather*}
\Psi(t,z)\simeq Q(\chi),\\
Q(\chi)=((\chi_1^2-\chi_2 ^2)+ \cdots+ (\chi_{2r-1}^2-\chi_{2r}^2) + (\chi_{2r+1}^2+\chi_{2r+2}^2)+\cdots +(\chi_{2n-1}^2+\chi_{2n}^2)).
\end{gather*}
Applying the stationary phase method for the $\chi$ variables, which is legal since the remaining integration w.r.t. $t$ is of compact support, we get:
\begin{equation*}
\Upsilon_{z_0}(E_c,\varphi,h)\sim e^{\frac{i\pi}{4} \mathrm{sgn}(Q)} \int\limits_{\mathbb{R}} \hat{\varphi}(t) d\mu_t(z_0) +\mathcal{O}(h).
\end{equation*}
Here $\mathrm{sgn}(Q)=(n+n-r)-r=2(n-r)$, next by checking:
\begin{itemize}
\item The value of $\chi(t,z_0)$,
\item The value of the Jacobian $\frac{D \chi}{Dz}$ at the point $z_0$,
\end{itemize}
it comes out, $\chi$ being given by
$\chi(t,z)=(d\Phi_t(z_0)-\mathrm{Id}))z+\mathcal{O}(||z||^2)$ at
the first order, and as long as
$\mathrm{det}(d\Phi_t(z_0)-\mathrm{Id})\neq 0$, that:
\begin{equation}
d\mu_t(z_0)=| \det(d\Phi_t(z_0) -\mathrm{Id}) |^{-\frac{1}{2}}
\end{equation}
Observe that this coefficient, the Duistermaat-Guillemin-Uribe
density, is indeed a smooth function as long as we stay away from
any period of the linearized flow at $z_0$. An explicit
computation using $A(t)$ in our coordinates, done for example in
\cite{KhD1}, shows that the density is given by:
\begin{equation*}
d\mu_t(z_0)=\frac{1} {|\prod\limits_{j=1} ^r \mathrm{sinh}(\alpha_j(z_0) t) \prod\limits_{j=r+1}^{n} \sin(\alpha_j(z_0) t)|}.
\end{equation*}
On this formula we see that we have $r$-negative and $n-r$
positive eigenvalues at the critical point $x_0$. The desired
result follows since $d\mu_t(z_0)$ determines:
\begin{itemize}
\item The signature of the Hessian of $V$ at $x_0$.%
\item Eigenvalues $\alpha_j(x_0)$.
\end{itemize}
The second affirmation follows by Taylor expanding the density and evaluating it at different points.\medskip\\
\textbf{Comments.}\\
Such a density was first introduced by Duistermaat and Guillemin
\cite{D-G}. But it seems that the exploitation of this term to
describe contributions of critical points goes back to Guillemin
and Uribe \cite{G-U}. This kind of density can be extended to
Morse-Bott singularities for $V$ (see \cite{BPU,KhD1}). For a
strict minimum of the potential $x_0$ the shape of $d\mu_t(x_0,0)$
is fundamental to get better inverse spectral results. In
particular one has to check the presence of resonant coefficients.
See \cite{CdV-G} and \cite{G-U2} for improved results near a
minimum.

These densities can be continued as meromorphic-distributions with
singular support:
\begin{itemize}
\item At the origin, such a study is done in \cite{BPU}.
\item With singular at a period $T$ of $d\Phi_t(z_0)$, this is done in \cite{KhD1}.
\end{itemize}
For an operator which is not a Schr\"odinger operator some new
terms can generally appear at a period of $d\Phi_t(z_0)$ (see
\cite{Cam0}). All these facts strongly suggest that the 'pike
singularity' of $d\mu_t(z_0)$ as $t\rightarrow T$ near a period
$T$ of $d\Phi_t(z_0)$ should also describe the symbol. For
example, for a Schr\"odinger operator a double eigenvalue $w>0$ of
$d^2V(x_0)$ will generates a singularity of double magnitude in
$d\mu_t(z_0)$ at the point $T=\frac{2\pi}{w}$.
\medskip\\
\textbf{Example.} An important toy model is the case of an $n$-dimensional harmonic oscillator:
\begin{equation*}
P_h=-\frac{h^2}{2} \Delta + \frac{1}{2}\sum\limits_{k=1}^{n} w_k^2 x_k^2,\, w_k\neq 0\, \forall k.
\end{equation*}
This model is one of the few Hamiltonians that can be explicitly
solved. Then, for $t$ small\footnote{For $t$ large one can use
compositions and the stationary phase method. The result obtained
is then exact since the phase in quadratic w.r.t. $(x,y)$.} and
$t\notin \mathbb{Z}\frac{\pi}{w_k}$, the kernel of the propagator
$U_h(t)=\exp(\frac{it}{h} P_h)$ can be explicitly computed as:
\begin{gather*}
K(t,x,y)=\left( \prod\limits_{k=1}^n \frac{w_k}{2i\pi \sin(w_k t)}\right)^{\frac{1}{2}} e^{\frac{i}{h}S(t,x,y)},\\
S(t,x,y)=\sum\limits_{k=1}^n \frac{w_k}{\sin (w_k t)} \left( \frac{1}{2}\cos(w_k t)(x_k^2+y_k^2)-x_k y_k \right).
\end{gather*}
This clearly shows that the small $h$ behavior of $U_h(t)$
determines the eigenvalues $w_k$. This important example can be
perturbed to use Dyson-expansions and treat general potentials
near a minimum of the potential, see \cite{Hez}.
\subsection{Degenerate singularities.}
As seen in section 2 it is sufficient to study a micro-localized
problem:
\begin{equation*}
\Upsilon _{z_{0}}(E_{c},h,\varphi)=\frac{1}{2\pi }\mathrm{Tr}\int\limits_{\mathbb{R}}e^{i%
\frac{tE_{c}}{h}}\hat{\varphi}(t)\psi ^{w}(x,h D_{x})\mathrm{exp}(-\frac{it}{h}%
P_{h})\Theta (P_{h})dt.
\end{equation*}
Here $\psi \in C_{0}^{\infty }(T^{\ast }\mathbb{R}^{n})$ is
micro-locally supported near $z_{0}$ (cf section 2). If the
support of $\psi$ is chosen small enough it is relatively easy to
obtain a normal form for the phase function of the FIO
approximating $\mathrm{exp}(-\frac{it}{h}%
P_{h})\Theta (P_{h})$. For the convenience of the reader we recall
the contributions of equilibriums in the trace formula. We note
$\mathrm{S}(\mathbb{S}^{n-1})$ the surface of $\mathbb{S}^{n-1}$
and in the next two propositions it is understood that conditions
$(\mathcal{H}_1)$ to $(\mathcal{H}_3)$ are satisfied.
\begin{proposition}\label{minimum}
If $x_0$ is a local minimum the first new wave invariant attached to $x_0$ are given by:
\begin{equation*}
\Upsilon _{z_{0}}(E_{c},h,\varphi)\sim
h^{\frac{n}{2}+\frac{n}{2k}-n}
\sum\limits_{j,l\in\mathbb{N}^2} h^{\frac{j}{2}+\frac{l}{2k}}
\Lambda_{j,l}(\varphi ),
\end{equation*}
where the $\Lambda_{j,l}$ are some distributions. The first new wave-invariant, attached to the leading
coefficient is:
\begin{equation*}
h^{\frac{n}{2}+\frac{n}{2k}-n}
\frac{\mathrm{S}(\mathbb{S}^{n-1})}{(2\pi)^n}
\int\limits_{\mathbb{S}^{n-1}} |V_{2k}(\eta)|^{-\frac{n}{2k}}
d\eta \int\limits_{\mathbb{R}_{+} \times \mathbb{R}_{+}}
\varphi(u^2 +v^{2k}) u^{n-1} v^{n-1} dudv.
\end{equation*}
\end{proposition}
Observe the simplicity of the distributional-coefficient acting on
$\varphi$ obtained more or less by computing a volume. Up to the
spherical-mean, this coefficient is uniquely determined by $n$ and
$k$. For an unstable critical point the situation is more
complicated and we have:
\begin{proposition}\label{maximum1}
If $x_0$ is a local maximum we have :
\begin{equation*}
\Upsilon _{z_{0}}(E_{c},h,\varphi)\sim
h^{\frac{n}{2}+\frac{n}{2k}-n}
\sum\limits_{m=0,1}\sum\limits_{j,l\in\mathbb{N}^2}
h^{\frac{j}{2}+\frac{l}{2k}}\mathrm{log}(h)^m \Lambda
_{j,l,m}(\varphi ).
\end{equation*}
If $\frac{n(k+1)}{2k}\notin \mathbb{N}$, the first non-trivial new wave invariant is given by :
\begin{equation*}
h^{\frac{n}{2} +\frac{n}{2k}-n}%
\left\langle W_{n,k},\varphi\right\rangle%
\frac{\mathrm{S}(\mathbb{S}^{n-1})}{(2\pi)^n}
\int\limits_{\mathbb{S}^{n-1}}
|V_{2k}(\eta)|^{-\frac{n}{2k}}d\eta.
\end{equation*}
\end{proposition}
At first view the result seems to be the same, but the distributions $W_{n,k}$ of Proposition \ref{maximum1} are respectively given by :
\begin{gather*}
\left\langle W_{n,k},\varphi\right\rangle
=\int\limits_{\mathbb{R}}
(C_{n,k}^{+}|t|_{+}^{n\frac{k+1}{2k}-1}+C_{n,k}^{-}|t|_{-}^{n\frac{k+1}{2k}-1})
\varphi (t) dt,
\text{ if }n \text{ is odd},\\
\left\langle W_{n,k},\varphi\right\rangle
=C_{n,k}^{-}\int\limits_{\mathbb{R}} |t|_{-}^{n\frac{k+1}{2k}-1}
\varphi (t) dt,\text{ if }n \text{ is even}.
\end{gather*}
The other options are given by:
\begin{proposition}\label{maximum2}
If $\frac{n(k+1)}{2k}\in \mathbb{N}$ and $n$ is odd then the
top-order coefficients are:
\begin{equation*}
C_{n,k} \log (h)h^{\frac{n}{2}+\frac{n}{2k}-n}%
\frac{\mathrm{S}(\mathbb{S}^{n-1})}{(2\pi)^n}\int\limits_{\mathbb{S}^{n-1}}
|V_{2k}(\eta)|^{-\frac{n}{2k}} d\eta  \int\limits_{\mathbb{R}}
|t|^{n\frac{k+1}{2k}-1} \varphi (t) dt.
\end{equation*}
Finally, if $\frac{n(k+1)}{2k}\in \mathbb{N}$ and $n$ is even,
and we have :
\begin{equation*}
C_{n,k}^{\pm} h^{\frac{n}{2}+\frac{n}{2k}-n}%
\frac{1}{(2\pi)^n} \int\limits_{\mathbb{S}^{n-1}}
|V_{2k}(\eta)|^{-\frac{n}{2k}}d\eta \int\limits_{\mathbb{R}}
|t|^{n\frac{k+1}{2k}-1} \varphi (t) dt.
\end{equation*}
\end{proposition}
A careful examination of the proof shows that the last case in
Proposition \ref{maximum2} is similar to the first subcase of
Proposition \ref{maximum1} (with $n$ odd) since
$C_{n,k}^{+}=C_{n,k}^{-}$. But we refer to \cite{Cam4} for the
technical details.
\begin{remark} \rm{To emphasize the consistency of these results we mention that:
\begin{itemize}
\item $C_{n,k}$, $C_{n,k}^{\pm}$ are non-zero universal constants
depending only on $n$ and $k$. See \cite{Cam4} for an analytic formulation.
\item Such terms $h^\alpha$ and $h^\alpha \log(h)$, $\alpha\in \mathbb{Q}$ never appear if
$E$ is regular.
\end{itemize}}
\end{remark}
In this work we will mainly use the order w.r.t. $h$ of these
coefficients (and the constants appearing in the expansions). For
a detailed proof of Proposition \ref{minimum} see \cite{Cam3} and
for Propositions \ref{maximum1} and \ref{maximum2} see
\cite{Cam4}. The case $k=1$, i.e., quadratic singularities, can
also be retrieved from certain results of \cite{BPU} with some
support restrictions but, again, this is sufficient to attain our objectives.\medskip\\
\textbf{Asymptotic expansion at a regular energy level.} With
$(\mathcal{H}_5)$ and when the energy $E$ is regular, we have:
\begin{gather*}
\Upsilon(E,h,\varphi)%
\sim \frac{h^{1-n}}{(2\pi)^n} \mathrm{LVol}(\Sigma_E)
\hat{\varphi}(0)+\sum\limits_{j=1}^{\infty} h^{1-n+j}
c_j(\hat{\varphi})(0)\\
+\sum\limits_{\rho\in \Sigma_E} e^{\frac{i}{h}S_\rho} e^{i\pi
\mu_\rho/4} \sum\limits_{j=0}^{\infty}D_{\rho,j}
(\hat{\varphi})(T_\rho) h^j.
\end{gather*}
We refer to \cite{PU} for a proof. In the r.h.s. the sum concerns
periodic orbits $\rho$ of energy $E$ and is finite since
$\mathrm{supp}(\hat{\varphi})$ is compact. Here $S_\rho$,
$\mu_\rho$ and $T_\rho$ are resp. the action, the Maslov-index and
the period of the closed orbit $\rho$ and both $c_j$, $D_{\rho,j}$
are differential operators of order $j$. If $\varphi$ satisfies
$(\mathcal{H}_4)$ we have $c_j(\hat{\varphi})(0)=0$ for all
$j\in\mathbb{N}$ and for each regular value $E\in[E_1,E_2]$:
\begin{equation}\label{sum orbits}
\Upsilon(E,h,\varphi)%
\sim \sum\limits_{\rho\in\Sigma_E} e^{\frac{i}{h}S_\rho} e^{i\pi
\mu_\rho/4} \sum\limits_{j=0}^{\infty} D_{\rho,j}
(\hat{\varphi})(T_\rho) h^j.
\end{equation}
We accordingly obtain that this term is bounded and a fortiori:
\begin{equation}\label{trace non critical}
\Upsilon(E,h,\varphi)=\mathcal{O}(1),\text{ } \forall E\in
[E_1,E_2]\backslash\{E_c^1,...,E_c^l\}.
\end{equation}
This point will justify Corollary \ref{weak}.\medskip\\
Next, by Lemma \ref{periods}, we have $T_\rho\geq T$ uniformly
w.r.t. $E\in[E_1,E_2]$. Hence if $E$ is not critical and if in
addition $(\mathcal{H}_2)$ is satisfied the sum over the periods
of Eq.(\ref{sum orbits}) is simply 0 and in Eq.(\ref{trace non
critical}) we obtain in fact a bound $\mathcal{O}(h^\infty)$ (a fortiori $(\mathcal{H}_5)$ is not required in that situation).\medskip\\
\textbf{Asymptotics at a critical energy level.} For $E=E_c^m$
critical there is always a continuous contribution w.r.t. $t$ in
the spectral distribution showing up the presence of a new wave
invariant. A fortiori, a choice of $\hat{\varphi}$ flat at the
origin, or with a small compact support, does not erase this term.
We have:
\begin{equation*}
\Upsilon(E_{c}^m,h,\varphi)  \sim \sum\limits_{j=1}^{N_m}
f_j(h),
\end{equation*}
where $N_m$ is the number of equilibrium points on
$\Sigma_{E_{c}^m}$ and each $f_j(h)$ is given by the leading term
of Propositions \ref{minimum},\ref{maximum1} and \ref{maximum2}.$\hfill{\blacksquare}$\medskip\\
Note that the bottom of a symmetric double well (degenerate or not) gives a similar
answer as a single well of same nature. Hence without extra micro-local
considerations (e.g. suitably localized eigenfunctions estimates) it is difficult to distinguish
these 2 different settings.\medskip\\
\textbf{Proof of Theorem \ref{Main2}.} First, the micro-local
Weyl-law for regular energies:
\begin{equation*}
\Upsilon(E,h,\varphi) \sim (2\pi h)^{1-n} \hat{\varphi}(0)
\mathrm{Lvol}(\Sigma_E),
\end{equation*}
computes the dimension $n$. Now assume given a critical value
$E_c$ with a single critical point. The only choice of the
spectral function $\varphi$ allows to detect $E_c$ via the
singularity $f(h)$ of Theorem \ref{Main}. The knowledge of
$f(h)$ determines the order of the contribution. For example,
if :
\begin{equation*}
f(h)\sim C h^\alpha \log(h),
\end{equation*}
the critical point is a maximum and $\alpha$ computes the degree
$2k$ of the singularity. With $\hat{\varphi}$, the knowledge of
$k$ allows to compute the quantity :
\begin{equation*}
\int\limits_{\mathbb{R}} |t|^{n\frac{k+1}{2k}-1} \varphi (t) dt.
\end{equation*}
A fortiori $C$ determines the average of
$|V_{2k}|^{-\frac{n}{2k}}$ on $\mathbb{S}^{n-1}$. Without
$\log(h)$, the nature of the critical point can be detected by
a symmetry argument w.r.t. $\varphi$ since we a priori know $n$
and $k$. In view of Propositions \ref{minimum},\ref{maximum1},\ref{maximum2} we
can choose $\varphi$ odd, even, symmetric or non-symmetric w.r.t.
the origin to conclude. Note that if $\hat{\varphi}$ is not even
$\varphi$ is a priori complex valued. $\hfill{\blacksquare}$
\begin{remark}
\rm{The spherical average of $V_{2k}$ is a Jacobian in polar
coordinates around $x_0$. For example, by composition with
$e^{-|t|}$ we obtain:
\begin{equation*}
\int\limits_{\mathbb{R}^n} e^{-|V_{2k}(x)|}dx
=\frac{1}{2k}\Gamma(\frac{n}{2k}) \int\limits_{\mathbb{S}^{n-1}}
|V_{2k}(\eta)|^{-\frac{n}{2k}} d\eta.
\end{equation*}
The same result holds by integration of $f(V_{2k}(x))$, if $f\in
L^1(\mathbb{R}_{+},r^{\frac{n}{2k}-1}dr)$.}
\end{remark}
\begin{remark}
\rm{Enlarging the list of singularities would provide a bigger
"dictionary". The case of non-homogeneous singularities for $V$ is
still an open problem, in particular because the determination of
an explicit asymptotic expansion w.r.t. $h$ can be a difficult
analytic problem.}
\end{remark}
\begin{remark}
\rm{There are certainly many easy generalizations of Theorem
\ref{Main2} to integrable Schr\"odinger operators, e.g. when
$V(x_1,x_2)=V_1(x_1)+V_2(x_2)$.}
\end{remark}
\section{Extensions. Examples}
In this section we propose now several generalizations of the main results.
\subsection{Operators with sub-principal symbols.}
We will show, shortly, how to extend the result of Theorem
\ref{Main} to the case of an $h$-admissible operator. I mention several obvious motivations to problems involving operators with non-zero sub-principal symbols.\medskip\\
\textbf{Witten-Laplacians.} B. Helffer \& J. Sj\"ostrand for
Witten Laplacians, see e.g. \cite{Hel} for an overview and
references, have obtained recently many interesting results for
these operators. For example, the Witten-Laplacian on zero-forms
attached to the measure $e^{-f/h}$ is:
\begin{equation*}
\Delta_{f,h}^{(0)}=-h^2 \Delta+ \frac{1}{4} |\nabla
f(x)|^2 -\frac{h}{2} \Delta f(x), \text{ } f\in
C^{\infty}(\mathbb{R}^n),
\end{equation*}
whose symbol $p(x,\xi)=p_0(x,\xi)+h p_1(x,\xi)$ depends on $h$.\medskip\\
\textbf{Schr\"odinger operators on a manifold.} A Schr\"odinger
operator attached to a Laplace-Beltrami operator on a Riemannian
manifold $M$ and $h$-quantized by exterior multiplication:
\begin{equation*}
P_h=-h^2 \Delta_M +V(x),\, V\in C^\infty(M),
\end{equation*}
generally involves a sub-principal symbol. In local coordinates with a metric $G=g^{ij}$, $G^{-1}=g_{ij}$ and $g=\det{G}$ the operator is:
\begin{gather*}
-h^2 \sum\limits_{i,j} \sqrt{g} \frac{\partial}{\partial x_i} \frac{1}{\sqrt{g}}g_{ij} \frac{\partial}{\partial x_j}+V\\
= -h^2 \sum\limits_{i,j} g_{ij}(x) \frac{\partial}{\partial x_j}\frac{\partial}{\partial x_j}+V %
+h^2 \sum\limits_{i,j} \sqrt{g} \frac{\partial}{\partial x_i}(\frac{1}{\sqrt{g}}g_{ij})(x)\frac{\partial}{\partial x_j}.
\end{gather*}
Hence, in the sense of the $h$-calculus, we have $p_h:=p_0+h p_1$ with:
\begin{gather*}
p_0(x,\xi)=\sum\limits_{i,j} g_{ij}(x)\xi_j\xi_i+V (x),\\
p_1(x,\xi)= \sqrt{g}(x) \sum\limits_{i,j} \left(\frac{\partial}{\partial x_i} \frac{1}{\sqrt{g}}g_{ij}\right)_{|x} \, \xi_j.
\end{gather*}
Observe that $p_1=0$ at every point where $\xi=0$.\medskip\\
\textbf{General case.} More generally, it is possible to consider
$h$-admissible operators $P_h$ whose symbols are given by
asymptotic sums $p_h\sim \sum h^j p_j$ (e.g. interpreted as a
Borel sum w.r.t. $h$) with principal symbol
$p_0(x,\xi)=\xi^2+V(x)$ and a subprincipal symbol $p_1\neq 0$. Of
course in the formula for $p_0$ you can also replace $\xi^2$ by
the metric at $x$ if you want to do a
similar construction on a Riemannian-manifold (non-necessarily compact).\medskip\\
\textbf{Modification of the first transport equation.} In the
previous example the sub-principal symbol $p_1$ was non-zero. This
requires a light correction. Starting from the results of section
3 we proceed as follows. To each element $u_{h}$ of
$I(\mathbb{R}^{n},\Lambda )$ we can associate canonically a
principal symbol $e^{\frac{i}{h}S}\sigma
_{\mathrm{princ}}(u_{h})$, where $S$ is a function on $\Lambda $
such that $\xi dx=dS$ on $\Lambda$. In fact, if $u_{h}$ can
locally be represented by an oscillatory integral with amplitude
$a$ and phase $\varphi$, then we have $S=S_{\varphi }=\varphi
\circ i_{\varphi }^{-1}$ and $\sigma _{\mathrm{princ}}(u_{h})$ is
a section of $|\Lambda |^{\frac{1}{2}} \otimes M(\Lambda )$,
where:
\begin{itemize}
\item $M(\Lambda )$ is the
Maslov vector-bundle of $\Lambda$.%
\item $|\Lambda |^{\frac{1}{2}}$ is the bundle of half-densities on $\Lambda$.
\end{itemize}
When $p_{1}\neq 0$, in the global coordinates $(t,y,\eta )$ on $\Lambda$, the
half-density of $U_{h}(t)$ is given by :
\begin{equation}
\nu(t,y,\eta)=\exp (i\int\limits_{0}^{t}p_{1}(\Phi _{s}(y,-\eta ))ds)|dtdyd\eta |^{\frac{1%
}{2}}.\label{demi densite}
\end{equation}
For this expression, related to the resolution of the first
transport equation for the propagator, we refer to Duistermaat and
H\"{o}rmander \cite{D-H}. Accordingly, the F.I.O. approximating the
propagator has the amplitude :
\begin{equation*}
\tilde{a}(t,z)= a(t,z)\exp (i \int\limits_{0}^{t}
p_1(\Phi_s(z))ds).
\end{equation*}
Since $z_0$ is an equilibrium we have $p_1(\Phi_s(z_0))=p_1(z_0)$,
$\forall s$, and :
\begin{equation}
\tilde{a}(t,z_0)= \hat{\varphi}(t) e^{it p_1(z_0)}.
\end{equation}
\textbf{a)} If $p_1(z_0)=0$. This happens in many interesting
situations (in particular for a Laplace-Beltrami operator, see
above). Here the top order coefficients in the trace formula
remains the same, also for the new wave-invariants at a critical energy level.\\
\textbf{b)} If $p_1(z_0)\neq 0$. By Fourier inversion formula we
simply replace $\varphi(t)$ by $\varphi(t+p_1(z_0))$ in all
integral formulae of Propositions \ref{minimum},\ref{maximum1} and
\ref{maximum2}. Note that when using $(\mathcal{H}_4)$, this has
absolutely no effect for the mean values and hence on the
detection of the critical energy levels.
\subsection{Eigenfunction estimates approach.}
We inspect now the case of an energy surface supporting more than
one critical point. The method we use here is in reality much more
restrictive (physically and also from the point of view of
spectral theory) since it implicitly use eigenfunctions estimates
via a $C_0^\infty(T^*M)$-observable. These observable are bounded
operators on $L^2(M)$ (e.g. via a Calderon-Vaillancourt estimates)
and can be inserted in the trace, trace class
operators being an ideal.\medskip\\
Since everything below is local we can freely assume that
$M=\mathbb{R}^n$, if not we can use local coordinates given by the
exponential map. Let $K=p^{-1}(J)\subset T^{*}\mathbb{R}^n$ be
compact and:
\begin{equation*}
d_0=\frac{1}{2}\inf\limits_{i\neq j} d(z_i,z_j),
\end{equation*}
where $d$ is any distance on $T^*\mathbb{R}^n$. By construction,
each open ball $B(z,d_0)\subset T^{*}\mathbb{R}^n$ contains at
most 1 critical point for each $z\in K$. Clearly, we can cover a
compact neighborhood of $K$ by a finite number of balls
$B(z,d_0)$. With a partition of unity adapted to this covering we
obtain :
\begin{equation*}
\sum\limits_{j=1}^{N} \psi_j^w(x,h D_x) =\mathrm{Id}, \text{
on } C_0^{\infty}(K).
\end{equation*}
For each energy $E\in J$, we obtain:
\begin{equation*}
\mathrm{Tr}\int\limits_{\mathbb{R}} \hat{\varphi}(t)
\Theta(P_h)e^{\frac{i}{h}t(P_{h}-E)}dt =\sum\limits_{j=1}^{N}
\mathrm{Tr}\int\limits_{\mathbb{R}} \hat{\varphi}(t) \psi_j^w(x,h
D_x)\Theta(P_h)e^{\frac{i}{h}t(P_{h}-E)}dt.
\end{equation*}
Note that the r.h.s. is studied in section 2. By the same argument
as before, if $\Sigma_s\cap \mathrm{supp}(\psi_j)$ contains no
critical point we obtain :
\begin{equation*}
\mathrm{Tr}\int\limits_{\mathbb{R}} \hat{\varphi}(t)
\psi_j^w(x,h D_x)\Theta(P_h)e^{\frac{i}{h}t(P_{h}-E)}dt=%
\mathcal{O}(h^{\infty}).\\
\end{equation*}
And if there is exactly one critical point $z_0\in\Sigma_{E_c}$ in
$\mathrm{supp}(\psi_j)$ we have :
\begin{equation*}
\mathrm{Tr}\int\limits_{\mathbb{R}} \hat{\varphi}(t) \psi_j^w(x,h
D_x)\Theta(P_{h})e^{\frac{i}{h}t(P_{h}-E_c)}dt=
\psi_j(z_0){f}_j(h),
\end{equation*}
and by construction no cancellation can occur since $\psi_j(z_0)>0$.\\
\begin{remark}\rm{In Corollary \ref{weak} we have considered
$(\mathcal{H}_5)$ for the flow. A similar result holds for a chaotic
dynamics and an isolated degenerate closed orbit can be treated as
in \cite{Pop}. Finally, using the results of \cite{PU} one can
extend Corollary \ref{weak} to the case of families of periodic
orbits of dimension $d\leq n$.}
\end{remark}
It is important to notice that to put a pseudo-differential operator in the spectral estimates means that the inverse spectral problem
is now implicitly expressed in terms of some $L^2$-expectation:
\begin{gather*}
\mu^h_k(\psi_j)= \langle \varphi_k^h , \psi_j^w(x,hD_x)
\varphi_k^h \rangle,\\
P_h \varphi_k^h =\lambda_k(h) \varphi_k^h,\, \lambda_k(h)\in
[E-ch,E+ch],
\end{gather*}
attached to eigenvectors $\varphi_k^h$ of $P_h$, see, e.g., the
section 'eigenvector estimates' of \cite{BPU}. Recall that,
combining Egorov's theorem and Calderon-Vaillancourt estimates, in
the regime $h\rightarrow 0^+$ the measure $\mu^h_k$ becomes more
and invariant under the Hamiltonian flow. It is also sometimes
possible to obtain a measure concentrated on the critical-set (see
\cite{BPU}). There is no paradox here since the flow is constant
on the critical set of the principal symbol of $P_h$.

For a Schr\"odinger operator whose potential is not a Morse
function the inverse spectral problem seems to be ill-defined and
requires eigenvectors estimates (which are indeed much stronger
estimates than those based only the spectrum). Also in the
previous construction an interesting problem is to get an a priori
lower bound for the number $d_0$ without doing any iteration on
the size of $\mathrm{supp}(\psi_j)$.
\subsection{Pseudo-differential operators.}
We can also apply the previous strategy for an
$h$-pseudo-differential operator with an isolated homogenous
singularity as this was considered in \cite{Cam5}. In fact, we
will stick to the simpler case of a local extremum as considered
in \cite{Cam1}. Assume that $p\sim p_0+h p_1 +\mathcal{O}(h^2)$
and that $p_0$ has a unique critical point $z_{0}=(x_{0},\xi
_{0})$ on the critical energy surface $\Sigma_{E_c}$. Also near
$z_0$ we have a conical (homogeneous) singularity:
\begin{equation*}
p_0(z)=E_{c}+\sum\limits_{j=k}^{N}\mathfrak{p}_{j}(z)+\mathcal{O}(||(z-z_{0})||^{N+1}),\, k>2,
\end{equation*}
where the functions $\mathfrak{p}_{j}$ are homogeneous of degree
$j$ w.r.t. $z-z_{0}$. Furthermore, assume that $z_{0}$ is a local
extremum of $\mathfrak{p}_{k}$. This implies that the first
non-zero homogeneous component $\mathfrak{p}_{k}$ is even and is
positive or negative definite. A fortiori $z_0$ is isolated on
$\Sigma_{E_c}$. An elementary example in dimension 1 is:
\begin{equation*}
p(x,\xi)=\pm (\xi^4+x^4) +R(x,\xi),
\end{equation*}
where $R(x,\xi)=\mathcal{O}(||(x,\xi)||^5)$ and $R$ is chosen so
that $p$ is confining and with tempered growth. Applying the
results of \cite{Cam1} or \cite{Hel0} we can retrieve 2 invariants
of $p$:
\begin{prop} \label{Invariant Pseudo} The new wave invariants at $E=E_c$ determines:
\begin{equation}
A(p_0)=\frac{1}{(2\pi)^{n}}\int\limits_{\mathbb{S}^{2n-1}}|\mathfrak{p}_{k}(\theta
)|^{-\frac{2n}{k}}d\theta,
\end{equation}
and the degree of homogeneity $k$ of $\mathfrak{p}_{k}$.
\end{prop}
\noindent\textbf{Proof.} The method of proof is the same as
before, just select $\hat{\varphi}$ such that all the usual
wave-invariants:
\begin{itemize}
\item Energy-surface distributions,%
\item Distribution supported by periodic-orbits,%
\end{itemize}
disappear from the spectral estimates $\Upsilon(E,h,\varphi)$. For
the new-wave invariants, following the construction of \cite{Cam1}
or \cite{Hel0}, it is possible to transform locally the phase
$\Psi$ of our FIO into $-tp_k$. The new asymptotic problem has the
form:
\begin{equation*}
\Upsilon_{z_0}(E_c,h,\varphi)\sim \frac{1}{(2\pi
h)^n}\int\limits_{\mathbb{R}\times T^*\mathbb{R}^n}
e^{-\frac{it}{h} p_k(x,\xi)} w(t,x,\xi)dtdxd\xi, \, w\in
C_0^\infty.
\end{equation*}
After integration w.r.t. $t$ the becomes elliptic and a little
discussion concerning the asymptotic behavior of oscillatory
integrals  show that:
\begin{equation*}
\Upsilon_{z_0}(E_c,h,\varphi)= h^{\frac{2n}{k}-n} \Lambda
_{0,k}(\varphi ) + \mathcal{O}(h^{\frac{2n+1}{k}-n}).
\end{equation*}
This result determines the even integer $k$. It is possible to
give a full asymptotic expansion (in powers of $h^{\frac{1}{k}}$)
but the first new wave-invariant at $E=E_c$ is explicitly given by
the distribution:
\begin{equation*}
\Lambda _{0,k}(\varphi ) =\frac{1}{k} \left\langle\varphi
(t+p_{1}(z_{0})),t_{z_0}^{\frac{2n-k}{k}}\right\rangle
\frac{1}{(2\pi
)^{n}}\int\limits_{\mathbb{S}^{2n-1}}|\mathfrak{p}_{k}(\theta
)|^{-\frac{2n}{k}}d\theta,
\end{equation*}
with $t_{z_0}=\mathrm{max}(t,0)$ if $z_0$ is a minimum and
$t_{z_0}=\mathrm{max}(-t,0)$ for a maximum. The 1rst coefficient
depends on $\varphi$, the dimension $n$, the degree $k\in
2\mathbb{N}^*$ of the singularity and the nature of the extremum.
The 2nd coefficient, independent of $\varphi$, determines the
spherical mean of $p_k$. $\hfill{\blacksquare}$\medskip\\
Two different spectral-estimates, with too functions $\varphi_1$
and $\varphi_2$ determine the ratio:
\begin{equation*}
r(\varphi_1,\varphi_2)=\frac{\left\langle\varphi_1
(t+p_{1}(z_{0})),t_{z_0}^{\frac{2n-k}{k}}\right\rangle}
{\left\langle\varphi_2
(t+p_{1}(z_{0})),t_{z_0}^{\frac{2n-k}{k}}\right\rangle },
\end{equation*}
so that the spherical mean on the sphere can generally be
determined after 2 spectral estimates.

Note that many symbols would give the same value. In particular
symbols conjugated by a rotation around the critical point cannot
be distinguished from this new-wave invariant.\medskip\\
\textbf{An invariance under re-scaling.} Now when $z_0$ is a
minimum, we observe that for any $j\in\mathbb{N}^*$ we can
re-scale our operator via:
\begin{equation*}
Q_{h,j}= (P_h-E_c)^j,
\end{equation*}
so that $z_0$ is still a minimum of the principal symbol $q_0=(p_0-E_c)^j$ of $Q_{h,j}$. Since:
\begin{equation*}
q_0(z)=(\mathfrak{p}_{k})^j(z)+\mathcal{O}(||z-z_0||^{kj+1}),
\end{equation*}
by a new application of Proposition \ref{Invariant Pseudo}, at the critical value zero, we get:
\begin{equation*}
A(q_0)=\frac{1}{(2\pi)^{n}}\int\limits_{\mathbb{S}^{2n-1}}(\mathfrak{p}_{k}(\theta
)^j)^{-\frac{2n}{kj}}d\theta =A(p_0).
\end{equation*}
It follows that $A(p_0)$ is also invariant under re-scaling (this argument does not apply for a local maximum).\medskip\\
\textbf{Acknowledgments.} This work was partially supported by a
Deutsche Forschungsgemeinschaft (D.F.G., the German research
foundation) Grant 'Micro-local analysis and geometry'. The DFG is
greatly acknowledged
for this support. \medskip\\
As a final remark, I mention that there is a lot of information to
retrieve from the spectral data of certain particular
Schr\"odinger operators on compact surfaces or manifolds. Certain
potentials like smoothed geodesic distances (distance functions
are singular at conjugate points) or height functions might
describe nicely the manifold. We plan to investigate these kind of
inverse spectral problems in a future article.


\begin{thebibliography}{00}
\bibitem{A-M} R. Abraham and J.E. Marsden, \textit{Foundations of
mechanics}, second edition, Benjamin/Cummings Publishing Co.,
Inc., Advanced Book Program, Reading, Mass.(1978).
\bibitem{BB} R. Balian and C.Bloch, Solution of the Schr\"odinger equation in term
of classical paths, Annals of Physics \textbf{85} (1974) 514-545.
\bibitem{Ber-Shu} F.A. Berezin and M.A. Shubin, \textit{The Schr\"odinger Equation}, Mathematics and Its Applications \textbf{66},
Kluwer Academic Publishers.
\bibitem{BPU} R. Brummelhuis, T. Paul and A. Uribe, Spectral estimates arround a critical level,
Duke Mathematical Journal \textbf{78} (1995) no. 3, 477-530.
\bibitem{BU} R. Brummelhuis and A. Uribe, A semi-classical trace formula for Schr\"odinger
operators, Communications in Mathematical Physics \textbf{136}
(1991) no. 3, 567-584.
\bibitem{Cam0} B. Camus, A semi-classical trace formula at a non-degenerate critical level,
Journal of Functional Analysis \textbf{208} (2004), no. 2,
446-481.
\bibitem{Cam1} B. Camus, A semi-classical trace formula at a totally degenerate critical
level. Contributions of extremums, Communications in Mathematical
Physics \textbf{207} (2004) no. 2, 513-526.
\bibitem{Cam3} B. Camus, Semi-classical spectral estimates for Schr\"odinger operators at a
critical level. Case of a degenerate minimum of the potential.
Journal of Mathematical Analysis and Applications (2007).
\bibitem{Cam4} B. Camus, Semi-classical spectral estimates for Schr\"odinger operators at a
critical level. Case of a degenerate maximum of the potential. Journal of Differential Equations \textbf{226} (2006)
no. 1, 295-322.
\bibitem{Cam5} B. Camus, Spectral estimates for degenerate critical levels.
Journal of Fourier Analysis and Applications \textbf{12} (2006),
no. 5, 495-515.
\bibitem{Cdv} Y. Colin de Verdi\`ere, Spectrum of the Laplace operator and periodic geodesics: Thirty years after. Annales de l'institut Fourier
\textbf{57} (2007) no. 7, 2429-2463.
\bibitem{CdV-G} Y. Colin de Verdi\`ere and V. Guillemin, A semi-classical inverse problem I: Taylor expansions.
Geometric aspects of analysis and mechanics, 81-95, Progr. Math.
\textbf{292} (2011).
\bibitem{DUI1} J.J. Duistermaat, Oscillatory integrals Lagrange immersions and unfolding
of singularities, Communications on Pure and Applied Mathematics
\textbf{27} (1974) 207-281.
\bibitem{D-G} J.J. Duistermaat and V. Guillemin, The spectrum of
positive elliptic operators and periodic bicharacteristics,
Inventiones Mathematicae \textbf{29} (1975), 39-79.
\bibitem{D-H} J.J. Duistermaat and L. H\"ormander, Fourier Integral
Operators, Acta mathematica \textbf{128} (1972) no. 3-4, 183-269.
\bibitem{G-U} V. Guillemin and A. Uribe, Circular symmetry and the trace formula, Inventiones Mathematicae
\textbf{96} (1989), 385-423.
\bibitem{G-U2} V. Guillemin and A. Uribe, Some inverse spectral results for semi-classical Schr\"odinger operators.
\bibitem{GUT} M. Gutzwiller, Periodic orbits and classical quantization conditions,
J. Math. Phys. \textbf{12} (1971) 343-358.
\bibitem{Haa} F. Haake, Quantum signatures of chaos. With a foreword by H.Haken.
Second edition. Springer-Verlag, Berlin, (2001).
\bibitem{Hel0} B. Helffer, Th\'eorie spectrale pour des op\'erateurs globalement elliptiques.
Ast\'erisque \textbf{112}. Soci\'et\'e Math\'ematique de France
(1984).
\bibitem{Hel} B. Helffer, Semiclassical analysis, Witten
laplacians, and statistical mechanics. Series on Partial
Differential Equations and Applications, 1. World Scientific
Publishing Co., Inc., River Edge, NJ (2002).
\bibitem{Hej} D.A. Hejhal, The Selberg trace formula for
$PSL(2,\mathbb{R})$, Vol.1, Springer-Verlag, L.N.M. 548.
\bibitem{Hez} H. Hezari, Inverse spectral problems for Schr\"odinger operators,
Communications in Mathematical Physics \textbf{288} (2009), 1061-1088.
\bibitem{KhD1} D. Khuat-Duy, A semi-classical trace formula for Schr\"odinger
operators in the case of a critical energy level, Journal of
Functional Analysis \textbf{146} (1997) no. 2, 299-351.
\bibitem{Laz} V.F. Lazutkin, KAM theory and semiclassical approximations to eigenfunctions.
With an addendum by A.I.Shnirel'man. Results in Mathematics and
Related Areas (3) 24. Springer-Verlag, Berlin, (1993).
\bibitem{Lieb} E.H. Lieb, Lieb-Thirring Inequalities, Kluwer Encyclopedia of Mathematics,
Supplement vol.II, p. 311-313 (2000)
\bibitem{PU} T. Paul and A. Uribe, The semi-classical trace formula and propagation of wave packets,
Journal of Functional Analysis \textbf{132} (1995), no. 1,
192-249.
\bibitem{Pop} G. Popov, On the contribution of degenerate periodic
trajectories to the wave-trace. Communications in Mathematical
Physics \textbf{196} (1998), no. 2, 363-383.
\bibitem{[Rob]} D. Robert, Autour de l'approximation
semi-classique, Progress in mathematics Volume 68, Birkh\"{a}user
Boston, Inc., Boston, MA, (1987).
\bibitem{Sim2} B. Simon, Nonclassical eigenvalue asymptotics. Joural of Functional Analysis \textbf{53}
(1983), no. 1, 84-98.
\bibitem{Sj-Zw} J. Sj\"ostrand and M. Zworski, Quantum monodromy and semi-classical trace formulae.
Journal de math\'ematiques pures et appliqu\'ees \textbf{81}
(2002), 1-33.
\bibitem{Yor} J.A. Yorke, Periods of periodic solutions and the Lipschitz constant,
Proceedings of the American Mathematical Society \textbf{69}
(1969) 509-512.
\end{thebibliography}
\end{document}